\documentclass[a4paper,11pt]{amsart}
 \usepackage[english]{babel}
  \usepackage{amsmath,amsfonts,amssymb,amsthm,relsize,setspace,nicefrac, yhmath,amscd}
  \usepackage{amsrefs, mathrsfs, enumerate} 
 \usepackage[colorlinks, linkcolor=red, citecolor=blue, urlcolor=blue, hypertexnames=true]{hyperref}

  \usepackage{hyperref} 

\usepackage[active]{srcltx} 
 \usepackage{color,soul} 
\newtheorem{theorem}{Theorem}[section]

\newtheorem{lemma}[theorem]{Lemma}
\newtheorem{proposition}[theorem]{Proposition}

\theoremstyle{definition}

\newtheorem{remark}{Remark}

\newcommand{\be}{\begin{equation}}
\newcommand{\bel}[1]{\begin{equation}\label{#1}}
\newcommand{\ee}{\end{equation}}

\newcommand{\barr}{\begin{eqnarray}}
\newcommand{\earr}{\end{eqnarray}}
\newcommand{\bars}{\begin{eqnarray*}}
\newcommand{\ears}{\end{eqnarray*}}

\newtheorem{subn}{\name}

\newcommand{\bsn}[1]{\def\name{#1}\begin{subn}}
\newcommand{\esn}{\end{subn}}

\newtheorem{sub}{\name}[section]

\newcommand{\bs}{\begin{sub}}
\newcommand{\es}{\end{sub}}

\newcommand{\bth}[1]{\def\name{Theorem}
\begin{sub}\label{t:#1}}
\newcommand{\blemma}[1]{\def\name{Lemma}
\begin{sub}\label{l:#1}}
\newcommand{\bcor}[1]{\def\name{Corollary}
\begin{sub}\label{c:#1}}
\newcommand{\bdef}[1]{\def\name{Definition}
\begin{sub}\label{d:#1}}
\newcommand{\bprop}[1]{\def\name{Proposition}
\begin{sub}\label{p:#1}}

\newcommand{\BA}{\begin{array}}
\newcommand{\EA}{\end{array}}
\newcommand{\BAN}{\renewcommand{\arraystretch}{1.2}
\setlength{\arraycolsep}{2pt}\begin{array}}
\newcommand{\BAV}[2]{\renewcommand{\arraystretch}{#1}
\setlength{\arraycolsep}{#2}\begin{array}}
\newcommand{\BSA}{\begin{subarray}}
\newcommand{\ESA}{\end{subarray}}

\newcommand{\BAL}{\begin{aligned}}
\newcommand{\EAL}{\end{aligned}}
\newcommand{\BALG}{\begin{alignat}}
\newcommand{\EALG}{\end{alignat}}
\newcommand{\BALGN}{\begin{alignat*}}
\newcommand{\EALGN}{\end{alignat*}}
\def\angb<#1>{\langle #1 \rangle}

\def\N{\mathbb{N}}

\def\R{\mathbb{R}}

\let\.=\cdot
\let\0=\emptyset

\numberwithin{equation}{section}

\theoremstyle{definition}

\let\.=\cdot
\let\0=\emptyset

\newenvironment{formula}[1]{\begin{equation}\label{eq:#1}}
                       {\end{equation}\noindent}

\def\Fi#1{\begin{formula}{#1}}
\def\Ff{\end{formula}\noindent}

\setstcolor{red}

\begin{document}


\title[Radially Blow up/global solution of free boundary cooperative  system]{ON The characterization of BLOW-UP and GLOBAL RADIAL SOLUTION FOR FREE BOUNDARY SYSTEM WITH NONLINEAR INHOMOGENEOUS GRADIENT AND SOURCE TERMS}


\author{Hoang Huy Truong}
\address{Faculty of Mathematics and Applications, Saigon University, 273 An Duong Vuong st., Ward 3, Dist.5, Ho Chi Minh City, Viet Nam}
\email{thhuy@sgu.edu.vn}
\thanks{}

\author{Hoang-Hung Vo$^{*}$}
\address{Faculty of Mathematics and Applications, Saigon University, 273 An Duong Vuong st., Ward 3, Dist.5, Ho Chi Minh City, Viet Nam}
\email{vhhung@sgu.edu.vn}
\thanks{$^*$ Corresponding author}


\begin{abstract} In the recent decade, a lot of attention has been drawn to the question on  the global phenomena of free boundary problems. Aiming to understand the  blow-up and global phenomena incorporating with free boundary, in this paper, we obtain a complete characterization of blowup and global phenomena for radial solutions of   a two-free boundaries  system as follows : 
\begin{align}\nonumber
 \left\{\begin{array}{rl}
u_t(t,r)= \Delta u(t,r) -  \lambda(t,x)|\nabla u(t,r)|^{\alpha} + a(t,x)v^{p}(t,r),& t>0,\  0<r<h(t),\\
v_t(t,r) = \Delta v(t,r) - \lambda(t,x) |\nabla v(t,r)|^{\alpha}+ a(t,x)u^{p}(t,r), & t>0,\  0 < r < g(t),
\end{array}\right.
\end{align}
where $r = |x|,\ x \in \R^N$, $p, \alpha >1 $ are given constants and $\lambda(t,x), a(t,x)$ satisfy suitable prescribed growth conditions.  First, we  prove the existence, uniqueness and stability of local solution to the system.  Then, we  classify the blowup and global phenomena by establishing some relations between $\alpha$, $p$ and growth rate of the coefficients. In particular, if $1<\alpha < p$, we show that the global fast solution exists for the initial value sufficiently small and the global slow solution exists with a suitable initial value, while blow-up solutions hold for sufficiently large initial value. On the other hand, if $\alpha \geq p$ associated with a suitable comparison on  $\beta, p$ and $\alpha$, then there exist global solutions with nonnegative initial data of exponential decay. Here, the presence of inhomogeneous gradient absorption prevents  blow-up phenomena and thus the solutions are slanted to the global phenomena. Our approach is  being far different from the seminal works \cite{MF,PS,SBP,TZ, TVH1,TVH2}, where the authors can only handled the equations without gradient term or with constant coefficients by self-similar solution or energy method. To our knowledge, this is the first work revealing the influence of the inhomogeneous coefficients to the blow-up and global phenomena to the cooperative system with nonlinear gradient and different free boundaries.

\end{abstract} 





\maketitle

\tableofcontents

\section{Introduction and Main results}

The study of blow-up and global solutions for  reaction–diffusion equations is one of the central topic of partial differential equations.  In the past decades, much attentions have been laid on the equations or the systems of non-variational  form, especially the equations associated with a gradient terms 
\cite{MF,PS,PZ,AZ1,AZ2,ZX,SBP,TZ,TVH1, TVH2, RGP4} since it is not only more mathematically challenging and also opens many new fundamental and important problems in mathematics and physics. Inspired from the mentioned works, in this paper, we will study the nonnegative blow-up and global solutions $(u(t,r),v(t,r))$ to the system with two different free boundaries :
\begin{align}\label{eqs_main}
\left\{\begin{array}{l}
u_t = \Delta u - \lambda(t,x) |\nabla u|^{\alpha} + a(t,x)v^{p},\ \text{for}\ t>0,\  0<r<h(t),\\
v_t = \Delta v - \lambda(t,x) |\nabla v|^{\alpha}+ a(t,x)u^{p},\ \text{for}\ t>0,\  0 < r < g(t),\\
u_r(t,0) = v_r(t,0) = 0,\ \text{for}\ t >0,\\
u(t,r) =0,\ \text{for}\ r \geq h(t)\ \text{and}\ t>0,\ v(t,r)= 0,\ \text{for}\ r \geq g(t)\ \text{and}\ t>0, \\  
h^{\prime}(t) = -\mu u_r(t,h(t)),\ g^{\prime}(t)= -\eta v_r(t,g(t)),\ \text{for}\ t>0,\\
 h(0)= h_0>0,\ u(0,r) = u_0(r),\ \text{for}\ 0 \leq r \leq h_0,\ \text{and}\ u(0,r) = 0,\ \text{for}\  r > h_0,\\ 
 g(0)= g_0>0,\ v(0, r) = v_0(r),\ \text{for}\ 0 \leq r \leq g_0,\ \text{and}\ v(0,r) = 0,\ \text{for}\  r > g_0,
\end{array}\right.
\end{align}
where $r = |x|$, $x \in \R^N, N \geq 1$, $\Delta w = w_{rr} + \dfrac{N-1}{r} w_r$ is the usual Laplace operator acting on spherically symmetric functions and $r = h(t)$ and $r = g(t)$ are  different free boundaries to be determined later. Throughout this paper, we assume that $\alpha, p  > 1, h_0, g_0, \mu,\eta$ are given positive constants and the initial condition $(u_0(r),v_0(r))$ satisfies the following assumptions 
\begin{align}\label{assume_main} 
 \left\{\begin{array}{l}
 u_0 \in C^{2}([0,h_0]), v_0 \in C^{2}([0,g_0]),  u'_0(0) = v'_0(0) = 0,\\
 u_0(h_0) = 0, u_0(r) >0\ \text{and}\ u'_0(r) < 0, \ \text{for}\ r \in  (0,h_0], \\
v_0(g_0)=0, v_0(r) > 0 \ \text{and}\  v'_0(r) < 0,\ \text{for}\ r \in  (0,g_0].
\end{array}\right.
\end{align}
Moreover, we assume the following conditions on $\lambda(t,x)$ and $a(t,x)$ as follows :

\textbf{(H) Hypothesis}: The functions $\lambda (t,\cdot) , a(t,\cdot) \in C^{1}(\R^N)$ for all $t \geq 0$ and there exist $\gamma \leq 0$ and $\beta \leq 0$ such that
\begin{align}\label{assume_1} 
 \left\{\begin{array}{l}
(\textbf{i})\ \ a_1(|x|+1)^{\gamma} \leq a(t,x) \leq a_2(|x|+1)^{\gamma},\ \text{for}\ x \in \R^N,\ a_1, a_2 >0, \ \text{for all}\ t \geq 0.\\
\textbf{(ii})\ \ \lambda_1 (|x|+1)^{\beta} \leq \lambda(t,x) \leq \lambda_2 (|x|+1)^{\beta},\ \text{for}\ x \in \R^N, \lambda_1, \lambda_2 >0,\ \text{for all}\ t \geq 0.
\end{array}\right.
\end{align}

For the reader's convenience, we introduce some notions and notations that are used later. We denote by $T_{max}\equiv T_{max}(u_0, v_0) \in (0, \infty]$ the maximal existence time of the solution of problem  (\ref{eqs_main}), $h_{\infty}:= \lim\limits_{t \to \infty} h(t)$ and $g_{\infty}:= \lim\limits_{t \to \infty} g(t)$. If $T_{max} < \infty$, then the solution blows up in finite time $T_{max}$, in the sense that
\[\lim\limits_{t\to T_{\max}}\big(\|u(t,\cdot)\|_{L^{\infty}([0,h(t)])} +\|v(t,\cdot)\|_{L^{\infty}([0,g(t)])}\big) = \infty.\]
Otherwise, if $T_{max} = \infty$ and $h_\infty + g_\infty <\infty$, then we say that the solution is a global fast solution and the solution decays uniformly as $t \to \infty$ at an exponential rate, in addition, if $T_{max} = \infty$ and $h_\infty + g_\infty =\infty$, then the solution is called a global slow solution, and it decays as $t \to \infty$ at most polynomial rate. 


The pioneering work for this type of problem is due to Chipot and Weissler \cite{MF}, who used the energy method to investigate the blow up/global solutions for reaction-diffusion equation with nonlinear gradient absorption on fixed boundary domain 
\begin{align}\label{eqs_main_single_CW}
\left\{\begin{array}{ll}
u_t = \Delta u - b|\nabla u|^{q} + |u|^{p-1}u,\ & x\in \Omega, t>0\\
u = 0,\ & x\in \partial\Omega,t>0\\
u(x,0) = \phi(x) \geq 0, \ & x\in \Omega,
\end{array}\right.
\end{align}
where $b\in\{0,1\}$. Later, Souplet and Weissler \cite{PS} consider more general equation of (\ref{eqs_main_single_CW}) with $0<b\leq2$ by self-similar solution approach and describe precisely the asymptotic behavior of its radially symmetric profile. In fact, it was shown by Tersenov \cite{TA} that, in general, the blow-up phenomena is prevented by the presence of the convection term and  the precise asymptotic behaviour for global solution for equation with gradient absorption has been further studied in the nice works of R. G. Pinsky \cite{RGP1, RGP2, RGP4}.

 The blow-up solution for one-phase Stefan problem with a power-type reaction has first been studied by Ghidouche,  Souplet,  and  Tarzia  \cite{HPD}, their problem reads as follow
\begin{align}\label{eqs_main_single_SL}
\left\{\begin{array}{ll}
u_t - u_{xx} = u^{p}, &  t>0,\  0< x <s(t),\\
u_x(t,0) = u(t,x) =0, & t >0,\\
s^{\prime}(t) = - u_x(t,s(t)), & t>0,\\
 s(0)= s_0>0,\ u(0,x) = u_0(x), & 0 \leq x \leq s_0,
\end{array}\right.
\end{align}
where $p >1$. The authors impose an energy condition involving the initial data under which the solution blows up in infinite time in $L^\infty$ norm. Moreover, it was shown in \cite{HPD} that all global solutions are bounded and decay uniformly to 0, and that either: (i) the free boundary converges to a finite limit and the solution decays at an exponential rate, or (ii) the free boundary grows up to infinity and the decay rate is at most polynomial, which correspond to the  global fast and slow solutions. Later, Zhou and Lin \cite{PZ} investigated the blow-up solution of the Stefan problem with integral nonlocal reaction term and showed that the blowup occurs if the initial datum is sufficiently large, otherwise, the global  fast or slow solution exists depending on the $L^\infty$-norm of initial condition.

Recently, the blowup and global (fast and slow) solutions for the free boundary problem with a nonlinear absorption gradient in one dimensional space were completely characterized by $p$ and $q$ in the interesting work of Z. Zang and X. Zang \cite{ZX} 
 
\begin{align}\label{eqs_main_ZZ}       
\left\{\begin{array}{ll}
u_t - u_{xx} = u^{p} - \lambda |u_{x}|^{q}, & t>0,\  0< x <s(t),\\
u_x(t,0) = u(t,x) =0, & t >0,\\
s^{\prime}(t) = -\mu u_x(t,s(t)), & t>0,\\
 s(0)= s_0>0,\ u(0,x) = u_0(x), & 0 \leq x \leq s_0,
\end{array}\right.
\end{align}
where $p, q > 1$. The authors established the local well-posedness and comparison principle for (\ref{eqs_main_ZZ}). In addition, it was shown in \cite{ZX} that in the instance of $p > q > 1$, the solution blows up in finite time for the initial value sufficiently large; but if the initial data are small, there exist global fast solutions which decay exponentially. Moreover, the existence of global slow solution with a suitable initial value is proved. While in the instance of $q \geq p >1$, all solutions with nonnegative initial data of exponential decay exist globally.   

Next, it is worth mentioning  that R. G. Pinsky \cite{RGP1, RGP2} characterized the existence and nonexistence of global solution as well as the behavior of life span by a critical exponent $p^*$ depending on the growth rate of potential to  the following equation 
\begin{equation}\label{eqs_G.Pinsky}    
u_t = \Delta u + a(x)u^p\ \text{in}\ \R^{N}, t>0
\end{equation} 
under the assumption that $a(x)\gneqq 0$ is on the order $|x|^m$, for $ m \in (-2, \infty)$, or that $0 \lneqq a(x)\leq C|x|^{-2}$, $a \in C^{\alpha}(\R^N)$ and $ p > 1$. Beside this, in their celebrated work, D. Andreucci and E. Di Benedetto \cite{DE} focused on investigating  nonnegative global solutions of degenerate parabolic equations :
\begin{equation}\label{eps_DE}   
u_t - \Delta u^m = \dfrac{u^p}{(1+|x|)^\gamma}\quad\R^N \times (0, T), \textrm{ for some $T > 0$,}  
\end{equation}    
where the range of parameters is $m \geq 1, p >1, \gamma \in \R$, and the initial datum is assumed to be a nonnegative Radon measure satisfying a suitable  technical condition determined by the parameters $\alpha; m; N; p$. In \cite{DE}, the authors asserted the main $L^{\infty}_{loc}$-estimates and integral gradient estimates on the solution and trace their connection with the blow up time. In addition, the authors provided global existence as well as non existence results, discriminating the sub and super-critical cases. Furthermore, the qualitative properties of quasi-linear parabolic equations with nonlinear gradient terms have been investigated, namely, S. Lian et al \cite{SHCWX} studied the equation
\begin{equation}
        u_t = \text{div} (|Du|^{p-2}Du) + |u|^{q-1}u - \lambda |Du|^{l} \ \text{in} \ \R^{N} \times (0,T],   
\end{equation}
 with \ $p > 2, q \geq 1, l \geq 1$ and $\lambda$ are real constants. It is proved \cite{SHCWX} that when $1 \leq l \leq p-1$,  local-in-time weak solution exists in the sense of D. Andreucci and E. DiBenedetto \cite{DE}, while for $p/2 \leq l \leq p-1$ the existence and nonexistence of global (in time) solutions are obtained in various situations.  

In recent decades, the characterization of blow up and global phenomena to semilinear parabolic equations  has been intensively attracted by the community. It worth to pointing out that Giga et al. \cite{YG, YR} has first obtained a bound of global solutions as well as nondegeneracy of blow-up for semilinear heat equations.  Further development on this direction has been made in the well-quoted works \cite{RGP4, YG, YR,  PM, FH,NFJ}. Also, many works have been devoted to studying the blowup/global phenomena  and without gradient terms (see e.g. \cite{FH,RW,DMJ, RM, RW,  ALV, YKH, MM, MM1, TVH1}) characterized by critical exponents. The first work is due to Escobedo and  Herrero \cite{ MM, MM1}, who considered the semilinear  cooperative parabolic system
\begin{align}\label{eqs_main_system_EH}  
\left\{\begin{array}{ll}
u_t - \Delta u = v^{p},\ & t>0, x\in\R^N \\
v_t - \Delta v = u^{q}, & t>0, x\in\R^N
\end{array}\right.
\end{align} 
 where  $p,q \in(0,\infty)$, $u(0,x), v(0,x) \geq 0$. The authors proved that both functions $u(x, t)$ and $v(x, t)$ must blow up in finite time if 
 $pq >1$ and $(\max\{p, q\}+1)/(pq-1) \geq N/2$, while there exits global solution, if $0 < pq \leq 1$. Later, H. A. Levine \cite{ALV} studied the above system in a more general situation on $D \subseteq \R^N$, ($D$ was a cone or $D$ was the exterior of a bounded domain) with $pq >1$ and $p,q \geq 1$, where the Fujita-type global existence–global nonexistence theorem is proved. It was shown in \cite{ALV} that there exists a critical value $c(D)$ such that (i) if $(\max\{p,q\} + 1)/(pq - 1) \geq c(D)$, no nontrivial, nonnegative global solutions exist, while (ii) if $(\max\{p,q\} + 1)/(pq - 1) < c(D)$, both nontrivial global and nonglobal solutions exist. Moreover, when $D$ is $\R^N$ or a cone, (i) holds with equality, but when $D$ is the exterior of a bounded domain, the case of equality in (i) is open. An explicit formula for $c(D)$ is given in each case. The system (\ref{eqs_main_system_EH}) with free boundary conditions in one space dimension was considered by M. Wang and Y. Zhao \cite{MY}, where the comparison principle was used to study the blowup/global solution questions, namely, if $p > 1$, then this system exists globally in time, provided that initial functions $u_0$ and $v_0$ are suitably small, whereas the solution blows up in finite time in the case initial functions $u_0$ and $v_0$ are large enough. Otherwise, if $0 < p \leq 1 $, then there are global solution.  Castillo and M. Loayza \cite{RM} extended the problems in \cite{MM, MM1,ALV} to more general parabolic systems with the coefficients are time-dependent and homogeneous Dirichlet boundary condition, in which the authors found the conditions that guarantee the global existence and the blow-up in finite time of nonnegative solutions. In another development, Ghoul-Nguyen-Zaag \cite{TVH1} studied  the cauchy problem for a non-variational (including gradient term) semilinear parabolic system, in which the authors establish the existence of a corresponding solution $(u, v)$ to this system such that $u$ and $v$ blows up simultaneously at the finite time $T_{max}$ with $a$ being the only blowup point of $u$ and $v$. However, there are still very limited investigations on the blow up and global solutions for free boundary problem, especially for the cooperative system in higher dimensional space with nonlinear gradient, and that is also the main goal of this paper.

Inspired by the mentioned works, in this paper, we investigate the influence of the coefficients to the blowup and global phenomena of the radial solutions of the two-free boundaries system (\ref{eqs_main}). Moreover, the presence of nonlinear gradient terms in the system (\ref{eqs_main}) yields significant difficulties since it requires deep technique to look for the sub and super solutions.  Here,  intriguingly, we find a profound effect of the coefficients associated with power-type nonlinearities on the blowup and global phenomena, which can be understood that the solution of (\ref{eqs_main})  with the different boundaries and different sources starting by different initial data but lead to the same global phenomena or blow-up in a finite time.

 Our  main result reads as follows :

\begin{theorem}\label{theorem_main1}
Assume that \textbf{(H)} holds, $(u,v,h,g)$ is the positive solution of (\ref{eqs_main}) and $\alpha, p >1$. Let $s_0 = \min\{h_0, g_0\}$ and $(u_0(r),v_0(r)) = A(\phi(r),\varphi(r))$, where $A > 0$ and $\phi, \varphi$ satisfy (\ref{assume_main}).
\begin{enumerate}[(1)]
\item \textbf{Blow-up}\

  \hspace{0.31cm} If $\alpha < p$ and $A$ is sufficiently large, then $T_{max} < \infty$ and 
  \[\lim\limits_{t\to T_{\max}}\big(\|u(t,\cdot)\|_{L^{\infty}([0,h(t)])} +\|v(t,\cdot)\|_{L^{\infty}([0,g(t)])}\big) = \infty.\]
\hspace{0.6cm} Moreover, there exists $C >0$ depending on $s_0, \gamma, p, \phi$ and $\varphi$ such that 
$ T_{\max} \leq  C.A^{-(p-1)}.$ 
\item \textbf{Global solution}
	\begin{enumerate}[(i)]
	\item \emph{(Global fast solution)} If $\alpha < p$ and $A$ satisfies condition  
\begin{align}\label{condition_GFS}
 0 < A \leq \dfrac{3}{4} \min\{(\dfrac{1}{16a_2s_0^2})^{\frac{1}{p-1}}, \dfrac{h_0^2}{8ps_0^2\mu}, \dfrac{g_0^2}{8s_0^2\eta}\} (\phi(0) + \varphi(0))^{-1}, 
 \end{align}
 then $T_{\max} = \infty,$ $h_{\infty} + g_{\infty} < \infty$ and there exist $h, k >0$, $C_1 >0$ depending on $u_0$, $C_2>0$ depending on $v_0$ such that 
$$ \|u(t,\cdot)\|_{L^{\infty}([0, h(t)])} \leq C_1 e^{-kt}\ \text{and}\ \|v(t,\cdot)\|_{L^{\infty}([0, g(t)])} \leq C_2 e^{-ht}.$$
More precisely, $C_1 =\dfrac{4}{3} \|u_0\|_{\infty},\ C_2 = \dfrac{4}{3}\|v_0\|_{\infty},$ and $k = h = \dfrac{1}{16s_0^2}$.

	\item \emph{(Global slow solution)} If $\gamma = 0$ in \textbf{(H)} and $\alpha < p$ then there exists $A > 0$ such that ($u, v, h, g$) is a global slow solution with initial data $u_0 (r) = A\phi(r)$ and $v_0 (r) = A\varphi(r)$, i.e.  $T_{max} = \infty$ and  $h_{\infty}  + g_{\infty} = \infty.$ Moreover, there exists a positive constant $C$ such that
	\[\sup_{t\geq 0}\|u(t,\cdot)\|_{\infty}+\sup_{t\geq 0}\|v(t,\cdot)\|_{\infty} \leq C\] 
	where $C$ depends on $\|u_0\|_{C_1}, \|v_0\|_{C_1}$ and $h_0, g_0$.
	\item If $\alpha \geq p$ and $\gamma \leq \beta \dfrac{p-1}{\alpha-1}$, then $T_{max} = \infty$. 
 \end{enumerate}
\end{enumerate}
\end{theorem}

By this result, we first characterize  the global behaviors of solutions for the system (\ref{eqs_main})  by comparing the parameters as $\alpha, p, \beta$ and $\gamma$ as well as the initial data. When $1<\alpha < p$, both blow up  and global  solutions to (\ref{eqs_main}) may occur depending on the largeness of the initial data. It is worth to emphasizing that  the absorption gradient supports the global phenomena while the source terms support the blowing up at finite time. As $\alpha < p$, we show that the absorption gradient is dominated by source terms but the blowing up only holds with large initial data while the global fast solutions even happens if the initial data is small enough. On the other hand,  when $\alpha \geq p$ and the growth rate of the coefficients supporting blowing up is controlled by growth rate of the coefficients supporting global behavior in the sense that $\gamma \leq \beta (p-1)/(\alpha-1)$, then all solutions starting by nonnegative  exponential decay data must hold globally. Moreover, we can estimate the upper bound for $T_{max}$ with respect to the initial data, $T_{\max}(A\phi,A\varphi
) \leq C A^{-(p-1)}$, which may be called the lifespan for blow up solutions. The construction of sub/super solutions and a comparison principle for  system (\ref{eqs_main}) plays a central role in our analysis, which  necessarily use a deep idea by exploiting the Stampacchia truncation method.  

\textbf{Organization of the paper.} The paper is organized as follows. In Section \ref{section_preliminary}, the well-posedness of the solution of the system (\ref{eqs_main}): local existence and uniqueness of the solution (Proposition \ref{proposition_lc&uni}) are studied by using the contraction mapping theorem and then applying the $L^p$ arguments and the Sobolev embedding theorem to demonstrate the continuous dependence (Proposition \ref{proposition_cdoid}). We also prove  the comparison principle in this section. In Section \ref{section_Blowup}, we prove that  the blow up solutions hold in the case $1 < \alpha < p$.  In Section \ref{section_global}, we give the existence of global fast solutions  and global slow solutions occurs  when $p> \alpha >1$ and no blowup solution  holds if $\alpha \geq p > 1$.  
\section{Preliminary results}\label{section_preliminary}
\subsection{The Well-posedness}  
In this section, we first verify the local existence and  uniqueness of solution by straightening free boundary and applying the contraction mapping theorem and then we study continuous dependence on initial data of (\ref{eqs_main}).      
\begin{proposition}[Local existence and uniqueness]\label{proposition_lc&uni}
Suppose that $(u_0, v_0)$ satisfy (\ref{assume_main}), \textbf{(H)} holds and $\beta_0 \in (0,1)$, then there exists a positive number $T>0$ such that system (\ref{eqs_main}) admits a unique positive solution 
\[(u,v,h,g) \in C^{\frac{1+\beta_0}{2},1+\beta_0}(Q_{T_h})\times C^{\frac{1+\beta_0}{2},1+\beta_0}(Q_{T_g})\times C^{1+\frac{\beta_0}{2}}([0,T])\times C^{1+\frac{\beta_0}{2}}([0,T]),\]
and 
\begin{equation}\label{lc&uni_1}
\|u\|_{C^{\frac{1+\beta_0}{2},1+\beta_0}(Q_{T_h})}+\|v\|_{C^{\frac{1+\beta_0}{2},1+\beta_0}(Q_{T_g})}+\|h\|_{C^{1+\frac{\beta_0}{2}}([0,T])} + \|g\|_{C^{1+\frac{\beta_0}{2}}([0,T])} \leq C,
\end{equation}
with $Q_{T_h} =\{(t,r) \in \R^{2} : t \in [0,T], r \in [0,h(t)]\}, Q_{T_g} =\{(t,r) \in \R^{2} : t \in [0,T], r \in [0,g(t)]\}, C > 0$ and $T$ depend only on $\|u_0\|_{C^{2}([0,h_0])}, \|v_0\|_{C^{2}([0,g_0])},h_0, g_0$ and $\beta_0$.   
\end{proposition}
\begin{proof}
Motivated from \cite{ZX,PZ}, we first straighten the free boundary. Let $\xi_1(i)$ and $\xi_2(i)$ be two function in
$C^3[0,\infty)$ satisfies 
\begin{align*}
\begin{array}{lll}
\xi_1(i) =\left\{\begin{array}{ll}
1, & |i-h_0| < \frac{h_0}{4},\\
0, & |i - h_0| > \frac{h_0}{2},
\end{array}\right.
& \text{and} & 
\xi_2(i) =\left\{\begin{array}{ll}
1, & |i - g_0| < \frac{g_0}{4},\\
0, & |i - g_0| > \frac{g_0}{2}.
\end{array}\right.
\end{array}
\end{align*}
Moreover, $|\xi_1'(i)| < \frac{5}{h_0}$ and  $|\xi_2'(i)| < \frac{5}{g_0}$  for all $i$.
We perform the following transformations
$$(t,y_1 )\to (t,x_1) \ \text{where} \  x_1 = y_1+\xi_1(|y_1|)(h(t) - h_0) \dfrac{y_1}{|y_1|}, y_1 \in \R^N,$$ and 
$$(t,y_2 )\to (t,x_2) \ \text{where} \  x_2 = y_2+\xi_2(|y_2|)(g(t) - g_0) \dfrac{y_2}{|y_2|}, y_2 \in \R^N,$$
which induce the two transformations
$(t,i) \to (t,r) \ \text{with} \ r = i +\xi_1(i)(h(t) - h_0), i \geq 0$ for $0 < r < h(t)$ \text{and}  $(t,i) \to (t,r) \ \text{with} \ r = i +\xi_2(i)(g(t) - g_0), i \geq 0$ for $0 < r < g(t)$.
 
For fixed $ t \geq 0$ as long as  $|h(t) - h_0| \leq \dfrac{h_0}{8}$ ($|g(t) - g_0| \leq \dfrac{g_0}{8}$), the transformation $ x_1 \to y_1$ ($ x_2 \to y_2$) determined above is a diffeomorphism from $\R^N$ onto $\R^N$ and the induced transformation $i \to r$ for $0 < r < h(t)$ ($i \to r$ for $0 < r < g(t)$) is a diffeomorphism from $[0,\infty)$ onto $[0,\infty)$. Furthermore, it changes the free boundary $|x_1| = h(t)$ ($|x_2| = g(t)$) to the fixed sphere $|y_1| = h_0$ ($|y_2| = g_0$). 
\\
Direct calculations yield
for  $0 < r < h(t),$
\begin{align*}
\dfrac{\partial i}{\partial r} = \dfrac{1}{1
+\xi_1'(i) (h(t)-h_0)}:=\sqrt{A_1(h(t),i)},\\
-\dfrac{1}{h'(t)}\dfrac{\partial i}{\partial t} = \dfrac{\xi_1(i)}{1+\xi_1'(i) (h(t)-h_0)}:=C_1(h(t),i),
\\
\dfrac{\partial^2 i}{\partial r^2} = -\dfrac{\xi_1''(i)(h(t)-h_0)}{[1+\xi_1'(i) (h(t)-h_0)]^3}:=B_1(h(t),i), 
\end{align*}
 denote $\dfrac{(N-1)\sqrt{A_1}}{i+\xi_1(i)(h(t)-h_0)}:=D_1(h(t),i)$.\\
  Similarly, for $0 < r < g(t),$ we have 
\begin{align*}
\dfrac{\partial i}{\partial r} = \dfrac{1}{1
+\xi_2'(i) (g(t)-g_0)}:=\sqrt{A_2(g(t),i)},\\
 -\dfrac{1}{g'(t)}\dfrac{\partial i}{\partial t} = \dfrac{\xi_2(i)}{1+\xi_2'(i) (g(t)-g_0)}:=C_2(g(t),i),\\
 \dfrac{\partial^2 i}{\partial r^2} = -\dfrac{\xi_2''(i)(g(t)-g_0)}{[1+\xi_2'(i) (g(t)-g_0)]^3}:=B_2(g(t),i), 
\end{align*}
and denote $\dfrac{(N-1)\sqrt{A_2}}{i+\xi_2(i)(g(t)-g_0)}:=D_2(g(t),i)$.\\
Let us set

$ 
\begin{aligned}
&u(t,r) = u(t,i+\xi_1(i)(h(t)-h_0)):= \varphi(t,i),\ \lambda(t,x) = \lambda(t,y+\xi_1(|y|)(h(t) - h_0)\frac{y}{|y|}):= \lambda_{01}(t,y),\\
&a(t,x) = a(t,y+\xi_1(|y|)(h(t) - h_0)\frac{y}{|y|}):= a_{01}(t,y)\ \text{for}\ 0 < r < h(t),
\end{aligned}
$
\\
and 

$\begin{aligned}
&v(t,r) = v(t,i+\xi_2(i)(g(t)-g_0)):= \chi(t,i),\ \lambda(t,x) = \lambda(t,y+\xi_2(|y|)(g(t) - g_0)\frac{y}{|y|}):= \lambda_{02}(t,y),\\
& a(t,x) = a(t,y+\xi_2(|y|)(g(t) - g_0)\frac{y}{|y|}):= a_{02}(t,y)\ \text{for}\ 0 < r < g(t),
\end{aligned}$ 
\\
then
\begin{align*}
\begin{array}{lll}
u_t =\varphi_t+\varphi_i\dfrac{\partial i}{\partial t}=\varphi_t-C_1.h'(t).\varphi_i, 
& & v_t =\chi_t+\chi_i\dfrac{\partial i}{\partial t}=\chi_t-C_2.g'(t).\chi_i,\\
u_r=\varphi_i.\dfrac{\partial i}{\partial r} =\sqrt{A_1}.\varphi_i, & \text{and} &v_r=\chi_i.\dfrac{\partial i}{\partial r} =\sqrt{A_2}.\chi_i,\\
u_{rr}=\varphi_{ii}A_1+\varphi_i.B_1, & & v_{rr}=\chi_{ii}A_2+\chi_i.B_2,\\
|\nabla u|^{\alpha} = |u_r|^{\alpha} = \sqrt{A_1}^{\alpha}|\varphi_i|^{\alpha},& &|\nabla v|^{\alpha} = |v_r|^{\alpha} = \sqrt{A_2}^{\alpha}|\chi_i|^{\alpha},
\end{array}
\end{align*}
the free boundary system (\ref{eqs_main}) becomes 
 \begin{align}\label{lc&uni_2}
\left\{\begin{array}{l}
\varphi_t - A_1\varphi_{ii} -(B_1+h'C_1+D_1)\varphi_i = a_{01}\chi^{p} -\lambda_{01}\sqrt{A_1}^\alpha |\varphi_i|^{\alpha}\ \text{for}\ t>0,\  0<i<h_0,\\
\chi_t - A_2\chi_{ii} -(B_2+g'C_2+D_2)\chi_i = a_{02}\varphi^{p} -\lambda_{02}\sqrt{A_2}^\alpha |\chi_i|^{\alpha}\ \text{for}\ t>0,\  0<i<g_0,\\
\varphi_i(t,0) = \varphi(t,h_0)= 0, \chi_i(t,0) = \chi(t,g_0)= 0\ \text{for}\ t>0,\\
h^{\prime}(t) = -\mu\varphi_i(t,h_0); g^{\prime}(t) = -\eta \chi_i(t,g_0)\ \text{for}\ t>0,\\
\varphi(0,i) =\varphi_0(i)= u_0(i)\ \text{for}\ 0 \leq i \leq h_0;
\chi(0, i) = \chi_0(i) = v_0(i)\ \text{for}\ 0 \leq i \leq g_0.
\end{array}\right.
\end{align}
Obviously, $h^{*}= -\mu\varphi'_0(h_0) \geq 0; g^*=-\eta \chi'_0(g_0) \geq 0$ and for $0 < T \leq \min\{\dfrac{h_0}{8(1+h^{*})}, \dfrac{g_0}{8(1+g^{*})}\}$, we define
\begin{align*}
&Q_{1T}^0 = [0,T]\times[0,h_0], Q_{2T}^0 = [0,T]\times[0,g_0],\\
&\mathcal{O}_{1T} = \{\varphi \in C^{0,1}(Q_{1T}^0):\varphi \geq 0, \varphi(0,i) = u_0(i), \|\varphi - u_0\|_{C^{0,1}(Q_{1T}^0)} \leq 1\},\\ 
&\mathcal{O}_{2T} = \{\chi \in C^{0,1}(Q_{2T}^0):\chi \geq 0, \chi(0,i) = v_0(i), \|\chi - v_0\|_{C^{0,1}(Q_{2T}^0)} \leq 1\},\\
&\mathcal{O}_{3T} = \{h \in C^1([0,T]):h(0) = h_0,h'(0)=h^{*}, \|h' - h^*\|_{C([0,T])} \leq 1\}.\\
&\mathcal{O}_{4T} = \{g \in C^1([0,T]):g(0) = g_0,g'(0)=g^{*}, \|g' - g^*\|_{C([0,T])} \leq 1\}.
\end{align*}
Then, we have 
$W_T = \mathcal{O}_{1T} \times \mathcal{O}_{2T} \times \mathcal{O}_{3T} \times \mathcal{O}_{4T}$ is a closed convex set in $C^{0,1}(Q_{1T}^0) \times C^{0,1}(Q_{2T}^0) \times C^{1}([0,T])  \times C^{1}([0,T])$, and $W_T$ is a complete metric space with the metric 
\[d(X_1,X_2)=\|\varphi_1-\varphi_2\|_{C^{0,1}(Q_{1T}^0)}+\|\chi_1-\chi_2\|_{C^{0,1}(Q_{2T}^0)}+\|h'_1-h'_2\|_{C([0,T])}+\|g'_1-g'_2\|_{C([0,T])},\]
with $X_1=(\varphi_1,\chi_1,h_1,g_1),X_2=(\varphi_2,\chi_2,h_2,g_2) \in W_T$. \\
Firstly, it is easy to show that for arbitrary $(\varphi,\chi,h,g) \in W_T$, we have  
$|h(t) - h_0| = t|h'(\delta_1)|$ and $|g(t) - g_0| = t|g'(\delta_2)|$ with some $\delta_1, \delta_2 \in (0,t)$ by Lagrange's theorem. Moreover, $|h'(\delta_1)| \leq |h'(\delta_1) - h'(0)|+|h'(0)|$ and $|g'(\delta_2)| \leq |g'(\delta_2) - g'(0)|+|g'(0)|$ hence  
\begin{align*}
\begin{array}{ll}
|h(t) - h_0|  \leq T( \|h' - h^{*}\|_{C([0,T])}+|h^{*}|) \leq T(1+h^{*}) \leq \dfrac{h_0}{8},\\
|g(t) - g_0|  \leq T (\|g' - g^{*}\|_{C([0,T])}+|g^{*}| )\leq T(1+g^{*}) \leq \dfrac{g_0}{8}.
\end{array}
\end{align*}
Therefore, the transformations $(t,i) \to (t,r)$ and $ (t,j) \to (t,r)$ introduced at the beginning of the proof are well defined. For any $(\varphi,\chi, h, g) \in W_T$, we consider the following initial - boundary value problem
 \begin{align}\label{lc&uni_3} 
\left\{\begin{array}{ll}   
\overline{\varphi}_t - A_1\overline{\varphi}_{ii} -(B_1+h'C_1+D_1)\overline{\varphi}_i = a_{01}\chi^{p} -\lambda_{01}\sqrt{A_1}^\alpha |\varphi_i|^{\alpha},& t>0,\  0<i<h_0,\\    
\overline{\varphi}_i(t,0) = \overline{\varphi}(t,h_0)= 0, & t>0,\\   
\overline{\varphi}(0,i) =\varphi_0(i), & 0 \leq i \leq h_0,\\  
\end{array}\right.
\end{align}
 since $\varphi_0 \in W_{\kappa}^{2}((0,h_0))$, then applying $L^p$ theory and the Sobolev imbedding theorem introduced in \cite{LSV}, we can get that above problem (\ref{lc&uni_3}) admits a unique solution $\overline{\varphi} \in C^{\frac{1+\beta_0}{2},1+\beta_0}(Q_{1T}^0)$ with $\beta_0 = 1-(N+2)/\kappa$ and 
 \begin{equation}\label{lc&uni_4}
 \|\overline{\varphi}\|_{C^{\frac{1+\beta_0}{2},1+\beta_0}(Q_{1T}^0)} \leq C_1,
 \end{equation}
where $C_1 > 0$ is a constant depend on $h_0, \beta_0, \|u_0\|_{W_\kappa^2([0,h_0])}$ and $\|v_0\|_{W_\kappa^2([0,g_0])}$. \\

Similarly, for given $(\varphi,\chi, h,g) \in W_T$ and  $\chi_0 \in W_{\kappa}^{2}((0,g_0))$, the following initial-boundary value problem 
\begin{align}\label{lc&uni_5}
\left\{\begin{array}{ll}
\overline{\chi}_t - A_2\overline{\chi}_{ii} -(B_2+g'C_2+D_2)\overline{\chi}_i = a_{02}\varphi^{q} -\lambda_{02}\sqrt{A_2}^\alpha |\chi_i|^{\alpha},& t>0,\  0<i<g_0,\\
\overline{\chi}_i(t,0) = \overline{\chi}(t,g_0)= 0, & t>0,\\
\overline{\chi}(0,i) =\chi_0(i), & 0 \leq i \leq g_0,\\ 
\end{array}\right.
\end{align}
has a unique solution $\overline{\chi} \in C^{\frac{1+\beta_0}{2},1+\beta_0}(Q_{2T}^0)$ and 
\begin{equation}\label{lc&uni_6}
 \|\overline{\chi}\|_{C^{\frac{1+\beta_0}{2},1+\beta_0}(Q_{2T}^0)} \leq C_1. 
\end{equation}
Defining  
$\overline{h}(t)=h_0 - \mu \int_{0}^{t}(\overline{\varphi}_i(\tau,h_0), \overline{g}(t)=g_0 -\eta\int_0^{t}\overline{\chi}_i(\tau,g_0))d\tau$ 
then $\overline{h}(0)=h_0,\ \overline{h}'(0) = h^*, \overline{g}(0)=g_0,$\ $\overline{g}'(0) = g^*,$ $\overline{h}'(t) = -\mu(\overline{\varphi}_i(t,h_0),$  and  $\overline{g}'(t)=- \eta\overline{\chi}_i(t,g_0)$, 
 \begin{equation}\label{lc&uni_6'}
  \|\overline{h}'\|_{C^{\beta_0/2}([0,T])} \leq \mu C_1 \leq C_2\ \text{and}\ \|\overline{g}'\|_{C^{\beta_0/2}([0,T])} \leq \eta C_1 \leq C_2.
 \end{equation}
Now, we introduce a mapping $F : W_T \rightarrow C^{0,1}(Q_{1T}^0) \times C^{0,1}(Q_{2T}^0) \times C^{1}([0,T])\times C^{1}([0,T])$ by 
$$F(\varphi,\chi,h,g) = (\overline{\varphi}, \overline{\chi},\overline{h},\overline{g}).$$

We claim that
\begin{align*}
&\|\overline{\varphi}-u_0\|_{C^{0,1}(Q_{1T}^0)} = \|\overline{\varphi}(t,i) - \overline{\varphi}(0,i)\|_{C^{0,1}(Q_{1T}^0)} \leq T^{\frac{1+\beta_0}{2}}\|\overline{\varphi}(t,i)\|_{C^{\frac{1+\beta_0}{2},1}(Q^0_{1T})},\\
&\|\overline{\chi}-v_0\|_{C^{0,1}(Q_{2T}^0)} = \|\overline{\varphi}(t,i) - \overline{\chi}(0,i)\|_{C^{0,1}(Q_{2T}^0)} \leq T^{\frac{1+\beta_0}{2}}\|\overline{\chi}(t,i)\|_{C^{\frac{1+\beta_0}{2},1}(Q^0_{2T})},\\
& |\overline{h}'(t)-\overline{h}'(0)| \leq \|h'\|_{C^{\frac{\beta}{2}}([0,T])}T^{\frac{\beta_0}{2}},\ |\overline{g}'(t)-\overline{g}'(0)| \leq \|g'\|_{C^{\frac{\beta_0}{2}}([0,T])}T^{\frac{\beta_0}{2}},
\end{align*}
thus, 
\begin{align*}
\max\{ \|\overline{\varphi}-u_0\|_{C(Q_{1T}^0)}, \|\overline{\chi}-v_0\|_{C(Q_{2T}^0)}\} \leq C_1T^{\frac{1+\beta_0}{2}};\ \ \ \max\{\|h'-h^{*}\|_{C([0,T])}, \|g'-g^{*}\|_{C([0,T])}\} \leq C_2T^{\frac{\beta_0}{2}}.
\end{align*}
So, if we choose $0 < T \leq T_0 := \min\big\{C_1^{\frac{-2}{1+\beta_0}}, C_2^{\frac{-2}{\beta_0}}, \dfrac{h_0}{8(1+h^{*})}, \dfrac{g_0}{8(1+g^{*})} \big\}$ then $F$ maps $W_T$ into itself $W_T$.  It will be show that $F$ is a contraction mapping on $W_T$ for sufficiently small $T > 0$.\\
Let $(\varphi_k, \chi_k, h_k, g_k) \in W_T$ for $k = 1, 2$ and denote $(\overline{\varphi}_k,\overline{\chi}_k, \overline{h}_k,\overline{g}_k ) = F (\varphi_k, \chi_k, h_k, g_k)$ for $k = 1, 2$. It is to obtain
$$ \| \overline{\varphi}_k\|_{C^{\frac{1+\beta_0}{2},1+\beta_0}(Q_{1T}^0)} \leq C_1, \ \|\overline{\chi}_k\|_{C^{\frac{1+\beta_0}{2},1+\beta_0}(Q_{2T}^0)} \leq C_1, \|\overline{h}_{k}\|_{C^{\frac{\beta_0}{2}}([0,T])} \leq C_2 \ \text{and} \ \|\overline{g}_{k}\|_{C^{\frac{\beta_0}{2}}([0,T])} \leq C_2,\ \text{for}\ k = 1,2.$$
Setting $\widehat{\varphi} = \overline{\varphi}_1 -\overline{\varphi}_2$ and $ \widehat{\chi} = \overline{\chi}_1 - \overline{\chi}_2$, then 
\begin{align}\label{lc&uni_7}
\left\{\begin{array}{lr}
\widehat{\varphi}_t - A_1(i,h_2)\widehat{\varphi}_{ii} -(B(i,h_2)+h_2'C(i,h_2)+D(i,h_2))\widehat{\varphi}_i = K(i,t),& t>0,\  0<i<h_0,\\
\widehat{\varphi}_i(t,0) = \widehat{\varphi}(t,h_0)= 0, & t>0,\\
\widehat{\varphi}(0,i) = 0, & 0 \leq i \leq h_0,\\ 
\end{array}\right.
\end{align}
where \begin{align*}
K (i,t)= & \ \ (A_1(i,h_1)-A_1(i,h_2))\overline{\varphi}_{1,ii}\\
&+[(B_1(i,h_1)-B_1(i,h_2))+(h'_{1}C_1(i,h_1)-h'_2C(i,h_2))+(D_1(i,h_1)-D_1(i,h_2))]\overline{\varphi}_{1,i}\\
&+a_{01}(\chi_1^{p} - \chi_2^{p}) +  \lambda_{01}\big[\sqrt{A_1(i,h_1)}^{\alpha} |\varphi_{1,i}|^{\alpha} - \sqrt{A_1(i,h_2)}^{\alpha}|\varphi_{2,i}|^{\alpha}\big]. 
\end{align*}  
By (\ref{lc&uni_4}), (\ref{lc&uni_6}), (\ref{lc&uni_6'}) and applying the $L^{p}$ estimates for parabolic equation and Sobolev's embedding theorem, it is deduced that 
\begin{equation}\label{lc&uni_8}
\|\widehat{\varphi}\|_{C^{\frac{1+\beta_0}{2},1+\beta_0}(Q_{1T}^0)} \leq C_{3}\big( \|\varphi_1 - \varphi_2\|_{C^{0,1}(Q_{1T}^0)}+\|\chi_1 - \chi_2\|_{C^{0,1}(Q_{2T}^0)}+\|h_1 - h_2\|_{C^1([0,T])} \big), 
\end{equation}
with $C_3 > 0$ which depends on $C_1, C_2$ and $A, B, C, D$ these are as defined above. Similarly, we also obtain  that
\begin{equation}\label{lc&uni_9}
\|\widehat{\chi}\|_{C^{\frac{1+\beta_0}{2},1+\beta_0}(Q_{2T}^0)} \leq C_{3}\big( \|\chi_1 - \chi_2\|_{C^{0,1}(Q_{2T}^0)}+\|\varphi_1 - \varphi_2\|_{C^{0,1}(Q_{1T}^0)}+\|g_1 - g_2\|_{C^1([0,T])} \big). 
\end{equation}
Moreover, 
\begin{align*}
\|\overline{h}'_1 - \overline{h}'_2\|_{C^{\frac{\beta_0}{2}}([0,T])} = \mu\| \widehat{\varphi}_i(\cdot,h_0)\|_{C^{\frac{\beta_0}{2},0}(Q_{1T}^0)} \leq \mu\| \widehat{\varphi}_i\|_{C^{\frac{\beta_0}{2},0}(Q_{1T}^0)} 
 \\
\|\overline{g}'_1 - \overline{g}'_2\|_{C^{\frac{\beta_0}{2}}([0,T])} = \eta\| \widehat{\chi}_i(\cdot,g_0)\|_{C^{\frac{\beta_0}{2},0}(Q_{2T}^0)} \leq \eta\| \widehat{\chi}_i\|_{C^{\frac{\beta_0}{2},0}(Q_{2T}^0)}. 
\end{align*}
As $\|\widehat{\varphi}_i\|_{C^{\frac{\beta_0}{2},0}(Q_{1T}^0)} \leq C_0\|\widehat{\varphi}\|_{C^{\frac{1+\beta_0}{2},1+\beta_0}(Q_{1T}^0)}$ and $\|\widehat{\chi}_i\|_{C^{\frac{\beta_0}{2},0}(Q_{2T}^0)} \leq C_0 \|\widehat{\chi}\|_{C^{\frac{1+\beta_0}{2},1+\beta_0}(Q_{2T}^0)}$, for some $C_0 >0$, we find
\[\|\overline{h}'_1 - \overline{h}'_2\|_{C^{\frac{\beta_0}{2}}([0,T])} +\|\overline{g}'_1 - \overline{g}'_2\|_{C^{\frac{\beta_0}{2}}([0,T])} \leq C_4 \big(\|\widehat{\varphi}\|_{C^{\frac{1+\beta_0}{2},1+\beta_0}(Q_T^0)}+ \|\widehat{\chi}\|_{C^{\frac{1+\beta_0}{2},1+\beta_0}(Q_T^0)}\big),\]
with $C_4 = (\mu +\eta)C_0$. We get that 
\begin{equation}\label{lc&uni_10}
\begin{array}{l}
\|\widehat{\varphi}\|_{C^{\frac{1+\beta_0}{2},1+\beta_0}(Q_{1T}^0)}+\|\widehat{\chi}\|_{C^{\frac{1+\beta_0}{2},1+\beta_0}(Q_{2T}^0)} + \|\overline{h}'_1 - \overline{h}'_2\|_{C^{\frac{\beta_0}{2}}([0,T])}+\|\overline{g}'_1 - \overline{g}'_2\|_{C^{\frac{\beta_0}{2}}([0,T])}\\
 \leq (1+C_{4})\big( \|\widehat{\varphi}\|_{C^{\frac{1+\beta_0}{2},1+\beta_0}(Q_T^0)}+ \|\widehat{\chi}\|_{C^{\frac{1+\beta_0}{2},1+\beta_0}(Q_T^0)}\big).
 \end{array}
\end{equation}
From (\ref{lc&uni_8}), (\ref{lc&uni_9}) and (\ref{lc&uni_10}), thus 
$$\begin{aligned}
& \|\widehat{\varphi}\|_{C^{\frac{1+\beta_0}{2},1+\beta_0}(Q_{1T}^0)}+\|\widehat{\chi}\|_{C^{\frac{1+\beta_0}{2},1+\beta_0}(Q_{2T}^0)} + \|\overline{h}'_1 - \overline{h}'_2\|_{C^{\frac{\beta_0}{2}}([0,T])}+\|\overline{g}'_1 - \overline{g}'_2\|_{C^{\frac{\beta_0}{2}}([0,T])}\\
& \leq (1+C_4)C_{3}\big( 2\|\varphi_1 - \varphi_2\|_{C^{0,1}(Q_{1T}^0)}+ 2\|\chi_1 - \chi_2\|_{C^{0,1}(Q_{2T}^0)}+\|h_1 - h_2\|_{C^1([0,T])}+\|g_1 - g_2\|_{C^1([0,T])}\big).
\end{aligned}$$
We see that $\|h_1 - h_2\|_{C^1([0,T])} \leq (T+1)\|h_1' - h_2'\|_{C([0,T])}$, $\|g_1 - g_2\|_{C^1([0,T])} \leq (T+1)\|g_1' - g_2'\|_{C([0,T])}$. Since $T \leq 1$, there holds 
$$\begin{aligned}
& \|\widehat{\varphi}\|_{C^{\frac{1+\beta_0}{2},1+\beta_0}(Q_{1T}^0)}+\|\widehat{\chi}\|_{C^{\frac{1+\beta_0}{2},1+\beta_0}(Q_{2T}^0)} + \|\overline{h}'_1 - \overline{h}'_2\|_{C^{\frac{\beta_0}{2}}([0,T])}+ \|\overline{g}'_1 - \overline{g}'_2\|_{C^{\frac{\beta_0}{2}}([0,T])}\\
& \leq (1+C_4)C_{3}\big( 2\|\varphi_1 - \varphi_2\|_{C^{0,1}(Q_{1T}^0)}+ 2\|\chi_1 - \chi_2\|_{C^{0,1}(Q_{2T}^0)}+(T+1)(\|h'_1 - h'_2\|_{C([0,T])}+\|g'_1 - g'_2\|_{C([0,T])})\big)\\
&\leq C_5\big( \|\varphi_1 - \varphi_2\|_{C^{0,1}(Q_T^0)}+ \|\chi_1 - \chi_2\|_{C^{0,1}(Q_T^0)}+\|h'_1 - h'_2\|_{C([0,T])}+\|g'_1 - g'_2\|_{C([0,T])}\big)\\
\end{aligned}$$
where $C_5 = 2(1+C_4)C_{3}$. 
Furthermore, 
\begin{align*}
&\|\overline{h}'_1 -\overline{h}'_2\|_{C([0,T])} \leq T^{\frac{\beta_0}{2}}\|\overline{h}'_1 -\overline{h}'_2\|_{C^{\frac{\beta_0}{2}}([0,T])}, \|\overline{g}'_1 -\overline{g}'_2\|_{C([0,T])} \leq T^{\frac{\beta_0}{2}}\|\overline{g}'_1 -\overline{g}'_2\|_{C^{\frac{\beta_0}{2}}([0,T])},\\
&\|\widehat{\varphi}\|_{C^{0,1}(Q_{1T}^0)} \leq T ^{\frac{1+\beta_0}{2}}\|\widehat{\varphi}\|_{C^{\frac{1+\beta_0}{2},1+\beta_0}(Q_{1T}^0)} \leq T^{\frac{\beta_0}{2}}\|\widehat{\varphi}\|_{C^{\frac{1+\beta_0}{2},1+\beta_0}(Q_{1T}^0)},\\
&\|\widehat{\chi}\|_{C^{0,1}(Q_{2T}^0)} \leq T^{\frac{1+\beta_0}{2}}\|\widehat{\chi}\|_{C^{\frac{1+\beta_0}{2},1+\beta_0}(Q_{2T}^0)} \leq T^{\frac{\beta_0}{2}}\|\widehat{\chi}\|_{C^{\frac{1+\beta_0}{2},1+\beta_0}(Q_{2T}^0)}. 
\end{align*}
It follows that 
\begin{align*}
&\|\widehat{\varphi}\|_{C^{0,1}(Q_{1T}^0)}+\|\widehat{\chi}\|_{C^{0,1}(Q_{2T}^0)}+ \|\overline{h}'_1 -\overline{h}'_2\|_{C([0,T])}+\|\overline{g}'_1 -\overline{g}'_2\|_{C([0,T])}\\
& \leq T^{\frac{\beta_0}{2}}\big(\|\widehat{\varphi}\|_{C^{\frac{1+\beta_0}{2},1+\beta_0}(Q_{1T}^0)}+ \|\widehat{\chi}\|_{C^{\frac{1+\beta_0}{2},1+\beta_0}(Q_{2T}^0)}+\|\overline{h}'_1 -\overline{h}'_2\|_{C^{\frac{\beta_0}{2}}([0,T])}+\|\overline{g}'_1 -\overline{g}'_2\|_{C^{\frac{\beta_0}{2}}([0,T])}\big)\\
& \leq T^{\frac{\beta_0}{2}}C_5\big( \|\varphi_1 - \varphi_2\|_{C^{0,1}(Q_{1T}^0)}+ \|\chi_1 - \chi_2\|_{C^{0,1}(Q_{2T}^0)}+\|h'_1 - h'_2\|_{C([0,T])}+\|g'_1 - g'_2\|_{C([0,T])}\big).
\end{align*}
Choosing $T = \min \{ 1,\big(\dfrac{1}{2C_5}\big)^{\frac{2}{\beta_0}}, (C_1)^{-\frac{2}{1+\beta_0}}, (C_2)^{-\frac{2}{\beta_0}},\dfrac{h_0}{8(1+h^{*})}, \dfrac{g_0}{8(1+g^{*})}\}$, we
see that the following estimates hold
\begin{align*}
&\|\overline{\varphi}_1-\overline{\varphi}_2\|_{C^{0,1}(Q_{1T}^0)}+\|\overline{\chi}_1-\overline{\chi}_2\|_{C^{0,1}(Q_{2T}^0)}+ \|\overline{h}'_1 -\overline{h}'_2\|_{C([0,T])}+\|\overline{g}'_1 -\overline{g}'_2\|_{C([0,T])}\\
& \leq \dfrac{1}{2}\big( \|\varphi_1 - \varphi_2\|_{C^{0,1}(Q_{1T}^0)}+ \|\chi_1 - \chi_2\|_{C^{0,1}(Q_{2T}^0)}+\|h'_1 - h'_2\|_{C([0,T])}+\|g'_1 - g'_2\|_{C([0,T])}\big).
\end{align*}  
With chosing $T$ is like above, we have $F$ is a contraction mapping on $W_T$, so $F$ has a unique fixed point $(\varphi, \chi, h, g) \in W_T$, which solves (\ref{lc&uni_2}), by the contraction mapping principle. In addition, we obtain that $\varphi \in C^{\frac{1+\beta_0}{2},1+\beta_0}(Q_{1T}^0), \chi \in C^{\frac{1+\beta_0}{2},1+\beta_0}(Q_{2T}^0), h, g \in C^{1+\frac{\beta_0}{2}}([0,T])$ and there exists a positive constant $C$ satisties 
$
\|\varphi\|_{C^{\frac{1+\beta_0}{2},1+\beta_0}(Q_{1T}^0)}+\|\chi\|_{C^{\frac{1+\beta_0}{2},1+\beta_0}(Q_{2T}^0)}+\|h\|_{C^{1+\frac{\beta_0}{2}}([0,T])} + \|g\|_{C^{1+\frac{\beta_0}{2}}([0,T])} \leq C.$ Furthermore, by the Schauder estimates, we have the higher regularity for $(u, v, h, g)$ as a solution of system (\ref{eqs_main}) such that $u \in C^{1+\frac{\beta_0}{2},2+\beta_0}(Q_{T_h})$, $v \in C^{1+\frac{\beta_0}{2},2+\beta_0}(Q_{T_g})$ and $h,g \in C^{1+\frac{\beta_0}{2}}([0,T])$, and (\ref{lc&uni_1}) holds. 
The proof of Proposition \ref{proposition_lc&uni} is complete.  
\end{proof}
\begin{proposition}[Continuous dependence on initial data]\label{proposition_cdoid} Assume that $(u_k, v_k,h_k,g_k,T_k)$ be the
unique solution of
\begin{align}\label{eqs_main_wp}
\left\{\begin{array}{l}
u_{k,t} = \Delta u_k - \lambda(t,x) |\nabla u_k|^{\alpha} + a(t,x)v_k^{p},\ \text{for}\ t>0,\  0<r<h_k(t),\\
v_{k,t} = \Delta v_k - \lambda(t,x) |\nabla v_k|^{\alpha}+ a(t,x)u_{k}^{p},\ \text{for}\ t>0,\  0 < r < g_k(t),\\
u_{k,r}(t,0) = v_{k,r}(t,0) = 0,\ u_k(t,h_k(t)) = v_k(t,g_k(t))= 0,\ \text{for}\ t>0, \\  
h_k^{\prime}(t) = -\mu u_{k,r}(t,h_k(t)),\ g_k^{\prime}(t)= -\eta v_{k,r}(t,g_k(t))),\ \text{for}\ t>0,\\
 h_k(0)= h_{0k}>0,\ u_k(0,r) = u_{0k}(r),\text{for}\ 0 \leq r \leq h_{0k},\\ 
 g_k(0)= g_{0k}>0,\ v_k(0, r) = v_{0k}(r),\text{for}\ 0 \leq r \leq g_{k0},
\end{array}\right.
\end{align}
with initial data $u_{0k}, v_{0k}, h_{0k}$ and $g_{0k}$ ($k = 1,2 $) and $T_k$ is the maximal existence time of the corresponding. Take $T =\min\{T_1, T_2\} \in (0, \infty]$. Then for all $0<T'<T$, there exists some constant $K>0$, depending on $\|u_{0k}\|_{C^2[0, m_0]}, \|v_{0k}\|_{C^2[0, n_0]}$, $\|u_k\|_{\mathcal{O}_{T1}}, \|v_k\|_{\mathcal{O}_{T2}}$, $h_{0k}, g_{0k}, \mu$ and $\eta$ such that 
\begin{center}
 $ \sup\limits_{t \in [0,T']} \|u_1(t,\cdot) - u_2(t,\cdot)\|_{C^0([0, m(t)])}+ \sup\limits_{t \in [0,T']} \|v_1(t,\cdot) - v_2(t,\cdot)\|_{C^0([0, n(t)])} + \|h_1 - h_2\|_{C^1([0,T'])}
+\|g_1 - g_2\|_{C^1([0,T'])}$\\
$ \leq K(\|u_{01}-u_{02}\|_{C^2([0,m_0])}+\|v_{01}-v_{02}\|_{C^2([0,n_0])}+|h_{01}-h_{02}|+|g_{01}-g_{02}|), $ 
 \end{center}
 where $m(t) = \min\{h_1(t),h_2(t)\}$, $n(t) = \min\{g_1(t),g_2(t)\}$, $m_0 =\min\{h_{01},h_{02}\}$, $n_0 =\min\{g_{01},g_{02}\}$, $\mathcal{O}_{T1} = [0,T]\times[0,m(t)]$ and $\mathcal{O}_{T2}=[0,T]\times[0,n(t)].$
\end{proposition}
\begin{proof} We first prove the local continuous dependence, i.e., there exists a $T_0>0$ and $K_0>0$
such that
\begin{equation}\label{wp_1}
\begin{aligned}
&\sup\limits_{t \in [0,T_0]} \|u_1(t,\cdot) - u_2(t,\cdot)\|_{C^0([0, m(t)])}+ \sup\limits_{t \in [0,T_0]} \|v_1(t,\cdot) - v_2(t,\cdot)\|_{C^0([0, n(t)])} + \|h_1 - h_2\|_{C^1([0,T_0])}+
\\
&+\|g_1 - g_2\|_{C^1([0,T_0])} \leq K(\|u_{01}-u_{02}\|_{C^2([0,m_0])}+\|v_{01}-v_{02}\|_{C^2([0,n_0])}+|h_{01}-h_{02}|+|g_{01}-g_{02}|).
\end{aligned}
\end{equation}  
We have $u_k \in C^{1+\frac{\beta_0}{2},2+\beta_0}([0, T_0]\times[0,h_k(t)]), v_k \in C^{1+\frac{\beta_0}{2},2+\beta_0}([0, T_0]\times[0,g_k(t)])$, $ h_k(t), g_k(t) \in C^{1+\frac{\beta_0}{2}}([0,T_0])$ by the existence and uniqueness of solutions to system (\ref{eqs_main}) and it implies that $h_k'(t) > 0$ and $g'_k(t) > 0$ by using Hopf's Lemma. Let us make the substitution $|y| = \dfrac{|x|}{h_k(t)}$ for $0 < |x| < h_k(t)$ and $|y| = \dfrac{|x|}{g_k(t)}$ for $0 < |x| < g_k(t)$. Also, we set $w_k(t,s) = u_k(t,h_k(t)s)$ and $z_k(t,s) = v_k(t,g_k(t)s)$, $k = 1, 2,$ with $s = |y|$.
 
Direct calculations yield
\begin{align}\label{wp_2}
\left\{\begin{array}{l}
w_{k,t} - \dfrac{\Delta w_{k}}{h_k^2} -\dfrac{h_k's}{h_k}w_{k,s} + \lambda(t,h_k s) \bigg|\dfrac{w_{k,s}}{h_k}\bigg|^{\alpha} = a(t,h_k s)z_k^{p},\ \text{for}\ 0 < t <T_0,\  0 < s  < 1,\\
z_{k,t} - \dfrac{\Delta z_{k}}{g_k^2} -\dfrac{g_k's}{g_k}w_{k,s} + \lambda(t,g_k s) \bigg|\dfrac{z_{k,s}}{g_k}\bigg|^{\alpha} = a(t,g_k s)w_k^{p},\ \text{for}\ 0 < t <  T_0,\  0 <s < 1,\\
w_{k,s}(t,0) = z_{k,s}(t,0) = 0,\ w_k(t,1) = z_k(t,1)= 0,\ \text{for}\ 0 < t <T_0, \\  
h_k^{\prime}(t) = -\mu \dfrac{w_{k,s}(t,1)}{h_k(t)},\ g_k^{\prime}(t)= -\eta \dfrac{z_{k,s}(t,1)}{g_k(t)},\ \text{for}\ 0 < t<T_0,\\
 h_k(0)= h_{0k}>0,\ w_k(0,s) = u_{0k}(h_{0k}s)= w_{0k}(s),\ \text{for}\ 0 \leq s \leq 1,\\ 
 g_k(0)= g_{0k}>0,\ z_k(0,s) = v_{0k}(h_{0k}s)= z_{0k}(s),\ \text{for}\ 0 \leq s \leq 1.
\end{array}\right.
\end{align}
Let us set $w(t,s) = w_1(t,s) - w_2(t,s)$ and $z(t,s) = z_1(t,s) - z_2(t,s)$  and then we obtain that 

 \begin{align}\label{wp_3}
\left\{\begin{array}{l}
w_{t} - \dfrac{\Delta w}{h_1^2} + A_1(t,s)w_{s}  = A_2(t,s)z + F_1(t,s),\ \text{for}\ 0 < t <T_0,\  0 < s  < 1,\\ 
w_{s}(t,0) =  w(t,1) = 0,\ \text{for}\ 0 < t <T_0,\\
w(0,s) = u_{01} - u_{02} = w_{0}(s),\ \text{for}\ 0 \leq s \leq 1,
\end{array}\right.
\end{align}
where 

$\begin{aligned}
&A_1(t,s) = \dfrac{\lambda(t, h_1s)\alpha\int_0^1\big|\xi(\frac{w_{1,s}}{h_1}-\frac{w_{2,s}}{h_2})+\frac{w_{2,s}}{h_2}\big|^{\alpha-2}\big[\xi(\frac{w_{1,s}}{h_1}-\frac{w_{2,s}}{h_2})+\frac{w_{2,s}}{h_2}\big]d\xi -h'_1s}{h_1},\\
&A_2(t,s) = a(t,h_1s)p\int\limits_{0}^{1}[\xi(z_1-z_2)+z_2]^{p-1}d\xi,\\
& F_1(t,s) = z_{2}^p(a(t,h_1s) - a(t,h_2s))+\bigg|\dfrac{w_{2,s}}{h_2}\bigg|^{\alpha}(\lambda(t,h_1s) - \lambda(t,h_2s))-\Delta w_{2}(\dfrac{1}{h^2_1}-\dfrac{1}{h^2_2})-s w_{2,s}(\dfrac{h'_1}{h_1}-\dfrac{h'_2}{h_2})
\\
& - \lambda(t,h_1s)\alpha w_{2,s}(\dfrac{1}{h_1}-\dfrac{1}{h_2})\int_0^1\bigg|\xi(\frac{w_{1,s}}{h_1}-\frac{w_{2,s}}{h_2})+\dfrac{w_{2,s}}{h_2}\bigg|^{\alpha-2}\bigg[\xi(\frac{w_{1,s}}{h_1}-\frac{w_{2,s}}{h_2})+\frac{w_{2,s}}{h_2}\bigg]d\xi,
\end{aligned}$

and  

 \begin{align}\label{wp_4}
\left\{\begin{array}{l}
z_{t} - \dfrac{\Delta z}{g_1^2} + B_1(t,s)z_{s}   =B_3(t,s)w + F_2(t,s),\ \text{for}\ 0 < t <T_0,\  0 < s  < 1,\\
z_{s}(t,0) =  z(t,1) = 0,\ \text{for}\ 0 < t <T_0, \\  
 z(0,s) = v_{01} - v_{02} = z_{0}(s),\ \text{for}\ 0 \leq s \leq 1,
\end{array}\right.
\end{align}
where 

$\begin{aligned}
&B_1(t,s) = \dfrac{\lambda(t, g_1s)\alpha\int_0^1\big|\xi(\frac{z_{1,s}}{g_1}-\frac{z_{2,s}}{g_2})+\frac{z_{2,s}}{g_2}\big|^{\alpha-2}\big[\xi(\frac{z_{1,s}}{g_1}-\frac{z_{2,s}}{g_2})+\frac{z_{2,s}}{g_2}\big]d\xi -g'_1s}{g_1},\\
&B_2(t,s) = a(t,g_1s)p\int\limits_{0}^{1}[\xi(z_1-z_2)+z_2]^{p-1}d\xi,\\
& F_2(t,s) = w_{2}^p(a(t,g_1s) - a(t,g_2s))+\bigg|\dfrac{z_{2,s}}{g_2}\bigg|^{\alpha}(\lambda(t,g_1s) - \lambda(t,g_2s))-\Delta z_{2}(\dfrac{1}{g^2_1}-\dfrac{1}{g^2_2})-s z_{2,s}(\dfrac{g'_1}{g_1}-\dfrac{g'_2}{g_2})\\
& - \lambda(t,g_1s)\alpha z_{2,s}(\dfrac{1}{g_1}-\dfrac{1}{g_2})\int_0^1\bigg|\xi(\frac{z_{1,s}}{g_1}-\frac{z_{2,s}}{g_2})+\dfrac{z_{2,s}}{g_2}\bigg|^{\alpha-2}\bigg[\xi(\frac{z_{1,s}}{g_1}-\frac{z_{2,s}}{g_2})+\frac{z_{2,s}}{g_2}\bigg]d\xi.
\end{aligned}$

Now, using $L^p$ arguments and the Sobolev's imbedding theorem, we can obtain
\begin{equation}\label{wp_5}
\begin{aligned}
\|w\|_{C^{\frac{1+\beta_0}{2},1+\beta_0}(D)} &\leq C_1(\|F_1\|_{L^p(D)} + \|z\|_{L^p(D)}+ \|w_0\|_{W^{2,0}_{p}([0,1])})\\
& \leq C_1(\|F_1\|_{C^0(D)} +\|z\|_{C^0(D)}+ \|w_0\|_{C^{2}([0,1])})\\
\end{aligned}
\end{equation}
and 
\begin{equation}\label{wp_6}
\begin{aligned}
\|z\|_{C^{\frac{1+\beta_0}{2},1+\beta_0}(D)} &\leq C_2(\|F_2\|_{L^p(D)} + \|w\|_{L^p(D)}+ \|z_0\|_{W^{2,0}_{p}([0,1])})\\
& \leq C_2(\|F_2\|_{C^0(D)}+ \|w\|_{C^0(D)} + \|z_0\|_{C^{2}([0,1])}),
\end{aligned}
\end{equation}
where $D = [0,T_0]\times[0,1]$, $C_1 >0$ and $C_2 > 0$ depend on $p, q, D$ and the modulus of continuity of $1/h^2_1$ and $1/g^2_1$.  
On the other hand, we see that 

\begin{equation}\label{wp_7}
\|F_1\|_{C^0(D)} \leq C_{01} \|h_1-h_2\|_{C^1([0,T_0])},
\end{equation}
and 
\begin{equation}\label{wp_8}
\|F_2\|_{C^0(D)} \leq C_{02} \|g_1-g_2\|_{C^1([0,T_0])},
\end{equation}
where $C_{01}>0$ and $C_{02} > 0$ depend on $\|u_{02}\|_{C^2([0,1])}$, $\|v_{02}\|_{C^2([0,1])}$, $h_{0i}$ and $g_{0i}$,  $i=1,2.$ \\
Moreover, 
\begin{equation}\label{wp_9}
\|z\|_{C^0(D)} \leq T_0^{\frac{\beta_0+1}{2}}\|z\|_{C^{\frac{1+\beta_0}{2},1+\beta_0}(D)} + \|z_0\|_{C^{2}([0,1])} 
\end{equation}
and 
\begin{equation}\label{wp_10}
\|w\|_{C^0(D)} \leq T_0^{\frac{\beta_0+1}{2}}\|w\|_{C^{\frac{1+\beta_0}{2},1+\beta_0}(D)} + \|w_0\|_{C^{2}([0,1])},
\end{equation}
we get that
\begin{equation}\label{wp_11}
\begin{aligned}
\|w_1 - w_2\|_{C^{\frac{1+\beta_0}{2},1+\beta_0}(D)}&\leq C_1(C_{01}\|h_1-h_2\|_{C^1([0,T_0])} +T_0^{\frac{\beta_0+1}{2}}\|z_1-z_2\|_{C^{\frac{1+\beta_0}{2},1+\beta_0}(D)}\\
&+\|z_{01}-z_{02}\|_{C^{2}([0,1])} + \|w_{01}-w_{02}\|_{C^{2}([0,1])})\\
\end{aligned}
\end{equation}
and 
\begin{equation}\label{wp_12}
\begin{aligned}
\|z_1 - z_2\|_{C^{\frac{1+\beta_0}{2},1+\beta_0}(D)}&\leq C_2(C_{02}\|g_1-g_2\|_{C^1([0,T_0])} +T_0^{\frac{\beta_0+1}{2}}\|w_1-w_2\|_{C^{\frac{1+\beta_0}{2},1+\beta_0}(D)} \\
&+\|z_{01}-z_{02}\|_{C^{2}([0,1])} + \|w_{01}-w_{02}\|_{C^{2}([0,1])}).\\
\end{aligned}
\end{equation}
We choose $0 < T_0 < \min\limits_{i=1,2}\{(\dfrac{1}{C_i})^{2/(\beta_0+1)}\}$, which implies that 

\begin{equation}\label{wp_13}
\begin{aligned}
&\|z_1 - z_2\|_{C^{\frac{1+\beta_0}{2},1+\beta_0}(D)} +\|w_1 - w_2\|_{C^{\frac{1+\beta_0}{2},1+\beta_0}(D)}\\
&\leq C_3\big(\|h_1-h_2\|_{C^1([0,T_0])}+ \|g_1-g_2\|_{C^1([0,T_0])} + \|z_{01}-z_{02}\|_{C^{2}([0,1])} + \|w_{01}-w_{02}\|_{C^{2}([0,1])}\big)\\
\end{aligned}
\end{equation}
thus, 
\begin{equation}\label{wp_14}
\begin{aligned}
&\sup\limits_{t \in [0,T_0]}\big(\|w_1(t,\cdot) - w_2(t,\cdot)\|_{C^{0}([0,1]} +\|z_1(t,\cdot) - z_2(t,\cdot)\|_{C^{0}([0,1])}\big) + \|h_1-h_2\|_{C^1([0,T_0])}+ \|g_1-g_2\|_{C^1([0,T_0])}\\
&\leq (C_3+1)(\|h_1-h_2\|_{C^1([0,T_0])}+ \|g_1-g_2\|_{C^1([0,T_0])}) + C_3(\|z_{01}-z_{02}\|_{C^{2}([0,1])} + \|w_{01}-w_{02}\|_{C^{2}([0,1])})\\
\end{aligned}
\end{equation}
where $C_3 > 0$ depends on $C_1, C_2$ and $\beta_0$. In addition, one has

\begin{equation}\label{wp_15}
\begin{aligned}
&\|h'_1-h'_2\|_{C^0([0,T_0])}+\|h_1-h_2\|_{C^0([0,T_0])} \leq (T_0+1)\|h'_1-h'_2\|_{C^0([0,T_0])}+ |h_{01} - h_{02}|,\\
&\|g'_1-g'_2\|_{C^0([0,T_0])}+\|g_1-g_2\|_{C^0([0,T_0])} \leq (T_0+1)\|g'_1-g'_2\|_{C^0([0,T_0])} + |g_{01} - g_{02}|,
\end{aligned}
\end{equation}
\begin{center}
$\begin{aligned}
&|h'_1 - h'_2| \leq \dfrac{\mu}{h_1}|w_{1,s}(t,1) - w_{2,s}(t,1)| + \dfrac{h'_2}{h_1}|h_1 - h_2|\\
&|g'_1 - g'_2| \leq \dfrac{\eta}{g_1}|z_{1,s}(t,1) - z_{2,s}(t,1)| + \dfrac{g'_2}{g_1}|g_1 - g_2|
\end{aligned}$
\end{center}
thus, 
\begin{equation}\label{wp_16}
\begin{aligned}
&\|h'_1 - h'_2\|_{C^0([0,T_0])} \leq \dfrac{2\mu}{h_{01}}T_0^{\beta_0/2}\|w_{1} - w_{2}\|_{C^{\frac{1+\beta_0}{2},1+\beta_0}(D)} + 2C_0|h_{01} - h_{02}|,\\
&\|g'_1 - g'_2\|_{C^0([0,T_0])} \leq \dfrac{2\eta}{g_{01}}T_0^{\beta_0/2}\|z_{1} - z_{2}\|_{C^{\frac{1+\beta_0}{2},1+\beta_0}(D)} + 2C_0|g_{01} - g_{02}|,
\end{aligned}
\end{equation}
since $T_0 < \min\{1,\dfrac{1}{2C_0}\} $, where $C_0 > 0$ depends on $h_{01}, g_{01}$ and $\beta_0$.  
From (\ref{wp_15}) - (\ref{wp_16}), we obtain that
\begin{equation}\label{wp_17}
\begin{aligned}
&\|h_1 - h_2\|_{C^1([0,T_0])} \leq \dfrac{2\mu}{h_{01}}T_0^{\beta_0/2}(T_0+1)\|w_{1} - w_{2}\|_{C^{\frac{1+\beta_0}{2},1+\beta_0}(D)} + (2C_0(T_0+1)+1)|h_{01} - h_{02}|,\\
&\|g_1 - g_2\|_{C^1([0,T_0])} \leq \dfrac{2\eta}{g_{01}}T_0^{\beta_0/2}(T_0+1)\|z_{1} - z_{2}\|_{C^{\frac{1+\beta_0}{2},1+\beta_0}(D)} + (2C_0(T_0+1)+1)|g_{01} - g_{02}|.
\end{aligned}
\end{equation} 

Setting $m =\max\{\dfrac{4\mu}{h_{01}}, \dfrac{4\eta}{g_{01}}\}$ and taking $T_0 = \min\{1, \dfrac{1}{2C_0}, (\dfrac{1}{C_1})^{2/\beta_0}, (\dfrac{1}{C_2})^{2/\beta_0}, \big(\dfrac{1}{2C_3(C_3+1)m}\big)^{2/\beta_0}\}$. 
\\
It follows from (\ref{wp_13}), (\ref{wp_14}) and (\ref{wp_17}) that

$
\begin{aligned}
&\sup\limits_{t \in [0,T_0]}\big(\|w_1(t,\cdot) - w_2(t,\cdot)\|_{C^{0}([0,1]} +\|z_1(t,\cdot) - z_2(t,\cdot)\|_{C^{0}([0,1])}\big) + \|h_1-h_2\|_{C^1([0,T_0])}+ \|g_1-g_2\|_{C^1([0,T_0])}\\
&\leq (C_3+1)(\dfrac{2\mu}{h_{01}}T_0^{\beta_0/2}(T_0+1)\|w_{1} - w_{2}\|_{C^{\frac{1+\beta_0}{2},1+\beta_0}(D)} + (2C_0(T_0+1)+1)|h_{01} - h_{02}|)\\
& \hspace{1cm} +(C_3+1)(\dfrac{2\eta}{g_{01}}T_0^{\beta_0/2}(T_0+1)\|z_{1} - z_{2}\|_{C^{\frac{1+\beta_0}{2},1+\beta_0}(D)} + (2C_0(T_0+1)+1)|g_{01} - g_{02}|)\\
&\hspace{1cm}+C_3(\|z_{01}-z_{02}\|_{C^{2}([0,1])} + \|w_{01}-w_{02}\|_{C^{2}([0,1])})\\
\end{aligned}
$

$
\begin{aligned}
&\leq (C_3+1)m T_0^{\beta_0/2}\big(\|w_{1} - w_{2}\|_{C^{\frac{1+\beta_0}{2},1+\beta_0}(D)}+\|z_{1} - z_{2}\|_{C^{\frac{1+\beta_0}{2},1+\beta_0}(D)}\big) \\
& + (C_3+1)(4C_0+1)(|g_{01} - g_{02}|+|h_{01} - h_{02}|)+C_3(\|z_{01}-z_{02}\|_{C^{2}([0,1])} + \|w_{01}-w_{02}\|_{C^{2}([0,1])})\\
\end{aligned}
$

$
\begin{aligned}
&\leq (C_3+1)m T_0^{\beta_0/2}C_3\big(\|h_1-h_2\|_{C^1([0,T_0])}+ \|g_1-g_2\|_{C^1([0,T_0])}\big)\\
&+(C_3+1)(4C_0+1)(|g_{01} - g_{02}|+|h_{01} - h_{02}|)\\
&+C_3((C_3+1)m T_0^{\beta_0/2}+1)(\|z_{01}-z_{02}\|_{C^{2}([0,1])} + \|w_{01}-w_{02}\|_{C^{2}([0,1])})\\
\end{aligned}
$

$
\begin{aligned}
&\leq \dfrac{1}{2}\big(\|h_1-h_2\|_{C^1([0,T_0])}+ \|g_1-g_2\|_{C^1([0,T_0])}\big)+(C_3+1)(4C_0+1)(|g_{01} - g_{02}|+|h_{01} - h_{02}|)\\
&+(C_3+\dfrac{1}{2})(\|z_{01}-z_{02}\|_{C^{2}([0,1])} + \|w_{01}-w_{02}\|_{C^{2}([0,1])}).
\end{aligned}
$

We deduce
\begin{equation}\label{wp_18}
\begin{aligned}
\sup\limits_{t \in [0,T_0]}&\big(\|w_1(t,\cdot) - w_2(t,\cdot)\|_{C^{0}([0,1]} +\|z_1(t,\cdot) - z_2(t,\cdot)\|_{C^{0}([0,1])}\big) + \|h_1-h_2\|_{C^1([0,T_0])}+ \|g_1-g_2\|_{C^1([0,T_0])}\\
&\leq K_1(\|w_{01}-w_{02}\|_{C^2([0,1])}+\|z_{01}-z_{02}\|_{C^2([0,1])}+|h_{01}-h_{02}|+|g_{01}-g_{02}|),
\end{aligned}
\end{equation} 
where $K_1 > 0$ depends on $C_0, C_1, C_2$ and $C_3$.
\\
Thus, if taking $T_0$ to be the initial time and
following the given proof, one finds $T_{01}>0$ and $K_{01}$ such that 

$
\begin{aligned}
&\sup\limits_{t \in [T_0,T_{01}]}\big(\|w_1(t,\cdot) - w_2(t,\cdot)\|_{C^{0}([0,1]} +\|z_1(t,\cdot) - z_2(t,\cdot)\|_{C^{0}([0,1])}\big) + \|h_1-h_2\|_{C^1([T_0,T_{01}])}+ \|g_1-g_2\|_{C^1([T_0,T_{01}])}\\
&\leq K_{01}(\|w_{1}(T_0,\cdot)-w_{2}(T_0,\cdot)\|_{C^2([0,1])}+\|z_{1}(T_0,\cdot)-z_{2}(T_0,\cdot)\|_{C^2([0,1])}+|h_{1}(T_0)-h_{2}(T_0)|+|g_{1}(T_0)-g_{2}(T_0)|).
\end{aligned}
$ 
Moreover, there must exist some constant $k >0$ such that 
\[|h_1(T_0)- h_2(T_0)|+|g_1(T_0)- g_2(T_0)| \leq k(|h_{01}-h_{02}|+|g_{01}-g_{02}|),\]
\[\|w_{1}(T_0,\cdot)-w_{2}(T_0,\cdot)\|_{C^2([0,1])}+\|z_{1}(T_0,\cdot)-z_{2}(T_0,\cdot)\|_{C^2([0,1])} \leq k (\|w_{01}-w_{02}\|_{C^2([0,1])}+\|z_{01}-z_{02}\|_{C^2([0,1])}).\]
Hence, we obtain that 
$$\begin{aligned}
&\sup\limits_{t \in [0,T_{01}]} \|w_1(t,\cdot) - w_2(t,\cdot)\|_{C^0([0, 1])}+ \sup\limits_{t \in [0,T_{01}]} \|z_1(t,\cdot) - z_2(t,\cdot)\|_{C^0([0, 1])} + \|h_1 - h_2\|_{C^1([0,T_{01}])}+
\\
&+\|g_1 - g_2\|_{C^1([0,T_{01}])} \leq K(\|w_{01}-w_{02}\|_{C^2([0,1])}+\|z_{01}-z_{02}\|_{C^2([0,1])}+|h_{01}-h_{02}|+|g_{01}-g_{02}|). 
\end{aligned}$$
By analogy, we get the conclusion
\end{proof}

\subsection{Comparison principle and Monotonicity}
In this section, we present the comparison principle which is very important to study the subsolution and supersolution of (\ref{eqs_main}). 
\begin{lemma}\label{lm_cp0}
Suppose that $T \in (0, \infty), h, g  \in C^{1}([0,T]), U \in C(\overline{\mathcal{D}^1_T})\cap C^{1,2}(\mathcal{D}^1_T)$, $V \in C(\overline{\mathcal{D}^2_T})\cap C^{1,2}(\mathcal{D}^2_T)$ are positive functions where $\mathcal{D}^1_T = \{ (t,r) \in \R^{2}: 0 < t \leq T, 0 < r < h(t)\}, \mathcal{D}^2_T = \{ (t,r) \in \R^{2}: 0 < t \leq T, 0 < r < g(t)\}$ and $(U, V, h, g)$ satisfies
\begin{align}\label{cp0_1}
\left\{\begin{array}{ll}
U_t - U_{rr} \leq  A_1|U_r|+ B_1U + C_1 V,\ \text{for}\ t>0,  0<r< h(t),\\
V_t - V_{rr} \leq  A_2|V_r|+ B_2V + C_2 U,\ \text{for}\ t>0, 0 < r < g(t),\\
U(t,r) = 0,\ \text{for}\ r \geq h(t)\ \text{and}\   V(t,r))= 0,\ \text{for}\ r \geq g(t), \\
 U_r(t,0) \geq 0,\ V_r(t,0) \geq 0,\ \text{for}\ t>0,\\  
h(0) = h_0, g(0) = g_0, U(0,r) \leq 0,\ \text{for}\ 0 \leq r \leq h_0;\ V(0,r) \leq 0,\ \text{for}\ 0 \leq r \leq g_0,
\end{array}\right.
\end{align}
where $A_i, |B_i|, C_i$ ($i =1,2$) are all nonnegative and bounded. Then $(U,V,h,g)$ satisfies for all $t \in [0,T]$,$$ U(t,r) \leq 0,\ \text{for}\ r \in [0,h(t)]\ \text{and}\  V(t,r) \leq 0,\ \text{for}\ r \in [0,g(t)].$$ 
\end{lemma}
\begin{proof} Let $\varepsilon > 0$, $\kappa = 2\ \max\limits_{i=1,2}\ \{ \sup_{D^i_T}B_i, \sup_{D^i_T}C_i\}$ and set $\mathcal{U} = U e^{-\kappa t} - \varepsilon
\psi$ and $\mathcal{V} = V e^{-\kappa t} - \varepsilon
\psi$, where
$$\psi (t,r) = M t + (1+r^2)^{1/2},$$
with $M = 1 + \max\limits_{i=1,2}\ \sup_{D^i_T}A_i$, $i=1,2$, then $\psi_{rr} - \psi_{t} +A_i |\psi_r|\leq 0$, $i =1,2$, it follows that a.e. in $D^1_T$, there holds
\begin{equation}\label{cp0_2}
\begin{aligned}
\mathcal{U}_{t} - \mathcal{U}_{rr}& = e^{-\kappa t}(U_t - U_{rr} - \kappa U) + \varepsilon(\psi_{rr} - \psi_{t}) \leq e^{-\kappa t}(A_1|U_r|+(B_1-\kappa)U + C_1 V)+ \varepsilon(\psi_{rr} - \psi_{t})\\
& \leq A_1|\mathcal{U}_r| + (B_1 - \kappa)\mathcal{U}+C_1\mathcal{V} +  \varepsilon(B_1 + C_1 - \kappa)\psi + \varepsilon(\psi_{rr}-\psi_t+A_1|\psi_r|)\\
& \leq A_1|\mathcal{U}_r| + (B_1 - \kappa)\mathcal{U}+C_1\mathcal{V}
\end{aligned}
\end{equation}
and similarly, a.e. in $D^2_T$,
\begin{equation}\label{cp0_3}
\mathcal{V}_{t} - \mathcal{V}_{rr} \leq A_2|\mathcal{V}_r| + (B_2 - \kappa)\mathcal{V}+C_2\mathcal{U}.
\end{equation} 
Setting $\mathcal{U}_{+} = \max\{0,\mathcal{U}\}$ and $\mathcal{V}_+=\max\{0,\mathcal{V}\}$, our assumptions imply $\mathcal{U}_{+} \in C(\mathcal{D}_{T}^1)\cap C^{1,2}(\mathcal{D}_{T}^1)$, $\mathcal{V}_{+} \in C(\mathcal{D}_{T}^2)\cap C^{1,2}(\mathcal{D}_{T}^2)$ $, \mathcal{U}_+(0,\cdot) = \mathcal{V}_+(0,\cdot) = 0 $ and for a.e. $t \in (0,T),\ \mathcal{U}_+(t)\in H_0^1(0,h), \mathcal{V}_+(t) \in H_0^1(0,g)$. Moreover, for a.e. $t \in (0,T)$, we have $\mathcal{U}_{rr}(\cdot,t)\in L^2(0,h),  \mathcal{V}_{rr}(\cdot,t)\in L^2(0,g) $, $(\mathcal{U}_+)_{r}(\cdot,t) = \chi_{\{\mathcal{U}>0\}} \mathcal{U}_{r}(\cdot,t),$ and $ (\mathcal{V}_+)_{r}(\cdot,t) = \chi_{\{\mathcal{V}>0\}}\mathcal{V}_{r}(\cdot,t)$. It follows from (\ref{cp0_2}) and $A_i, C_i$ ($i =1,2$) are nonnegative and bounded that

$ \begin{aligned}
&\dfrac{1}{2}\dfrac{d}{dt}\int_{0}^{h(t)}(\mathcal{U}_+)^2dr = \int_{0}^{h(t)}(\mathcal{U}_{+})_t \mathcal{U}_{+}dr + \dfrac{1}{2}h'(t)(\mathcal{U}_+)^2(t,h(t))\\
& \leq -\int_{0}^{h(t)}|(\mathcal{U}_{+})_r|^2 dr + \mathcal{P}·\int_{0}^{h(t)}|(\mathcal{U}_{+})_r| \mathcal{U}_{+}dr + \int_{0}^{h(t)}(B_1-\kappa)(\mathcal{U}_{+})^2dr + \mathcal{Q}\int_{0}^{h(t)} (\mathcal{V}_{+}) (\mathcal{U}_{+})dr
\end{aligned}$

$\begin{aligned}
& \leq -\int_{0}^{h(t)}|(\mathcal{U}_{+})_r|^2 dr + ·\int_{0}^{h(t)}| (\mathcal{U}_{+})_r|^2 dr + \dfrac{\mathcal{P}^2}{4}·\int_{0}^{h(t)} (\mathcal{U}_{+})^2dr+ \dfrac{\mathcal{Q}}{2}(\int_{0}^{h(t)} (\mathcal{U}_{+})^2dr + \int_{0}^{g(t)} (\mathcal{V}_{+})^2dr)
\end{aligned}
$

thus, 
\begin{equation}\label{cp0_4}
\dfrac{1}{2}\dfrac{d}{dt}\int_{0}^{h(t)}(\mathcal{U}_+)^2dr \leq (\dfrac{\mathcal{P}^2}{4}+\dfrac{\mathcal{Q}}{2})·\int_{0}^{h(t)} (\mathcal{U}_{+})^2dr+ \dfrac{\mathcal{Q}}{2}\int_{0}^{g(t)} (\mathcal{V}_{+})^2dr.
\end{equation}
Arguing as above, it follows from (\ref{cp0_3}) that
\begin{equation}\label{cp0_5}
\dfrac{1}{2}\dfrac{d}{dt}\int_{0}^{g(t)}(\mathcal{V}_+)^2dr \leq (\dfrac{\mathcal{P}^2}{4}+\dfrac{\mathcal{Q}}{2})·\int_{0}^{g(t)} (\mathcal{V}_{+})^2dr+ \dfrac{\mathcal{Q}}{2} \int_{0}^{h(t)} (\mathcal{U}_{+})^2dr,
\end{equation}
where $\mathcal{P}$ and $\mathcal{Q}$ are positive real numbers. 

Adding up (\ref{cp0_4}) and (\ref{cp0_5}), integrating, and using $\mathcal{U}_+(0,\cdot) =  \mathcal{V}_+(0,\cdot) = 0 $, we obtain that $\mathcal{U}_{+} = 0,$ in $\mathcal{D}^1_{T}$ and $ \mathcal{V}_{+} = 0$ in $\mathcal{D}^2_{T}$ and the conclusion follows by letting $\varepsilon \to 0.$
\end{proof}

\begin{lemma}\label{lm_cp}
Suppose that the condition \textbf{(H)} holds and $T \in (0, \infty), \overline{h}(t), \overline{g}(t)  \in C^{1}([0,T]), \overline{u} \in C(\overline{\mathcal{D}^1_T})\cap C^{1,2}(\mathcal{D}^1_T)$, $\overline{v} \in C(\overline{\mathcal{D}^2_T})\cap C^{1,2}(\mathcal{D}^2_T)$ are positive functions where $\mathcal{D}^1_T = \{ (t,r) \in \R^{2}: 0 < t \leq T, 0 < r < \overline{h}(t)\}, \mathcal{D}^2_T = \{ (t,r) \in \R^{2}: 0 < t \leq T, 0 < r < \overline{g}(t)\}$ and $(\overline{u}, \overline{v}, \overline{h}, \overline{g})$ satisfies
\begin{align}\label{cp_1}
\left\{\begin{array}{ll}
\overline{u}_t - \Delta \overline{u} + \lambda |\nabla \overline{u}|^{\alpha} \geq a\overline{v}^{p},\ \text{for}\ t>0,  0<r<\overline{h}(t),\\
\overline{v}_t - \Delta \overline{v} + \lambda |\nabla \overline{v}|^{\alpha} \geq a\overline{u}^{p},\ \text{for}\ t>0, 0 < r <\overline{g}(t),\\
 \overline{u}(t,r)= 0,\ \text{for}\ r \geq \overline{h}(t) ;\ \overline{v}(t,r)= 0,\ \text{for}\ r \geq \overline{g}(t),\\
\overline{u}_r(t,0) \leq 0,\ \overline{v}_r(t,0) \leq 0,\ \text{for}\ t>0, \\  
\overline{h}^{\prime}(t) \geq -\mu \overline{u}_r(t,\overline{h}(t)),\ \overline{g}^{\prime}(t) \geq -\eta \overline{v}_r(t,\overline{g}(t)),\ \text{for}\ t>0,\\
\overline{h}(0) = \overline{h}_0,\ \overline{g}(0) = \overline{g}_0, \overline{u}(0,r) = \overline{u}_0(r),\ \text{for}\ 0 \leq r \leq \overline{h}_0;\ \overline{v}(0,r) = \overline{v}_0(r),\ \text{for}\ 0 \leq r \leq \overline{g}_0. 
\end{array}\right.
\end{align}
If $ h_0 \leq \overline{s}_0,\ g_0 \leq \overline{g}_0, u_0(r) \leq \overline{u}_0(r)$ in $[0, h_0]\ \text{and}\ v_0(r) \leq \overline{v}_0(r)$ in $[0, g_0],$ 
then the solution $(u,v,h,g)$ of (\ref{eqs_main}) satisfies for all $t \in (0,T]$, 
\[h(t) \leq \overline{h}(t), g(t) \leq \overline{g}(t), u(t,r) \leq \overline{u}(t,r), r \in [0,h(t)]\ \text{and}\  v(t,r) \leq \overline{v}(t,r), r \in [0,g(t)].\] 
\end{lemma}
\begin{proof}
Motivated from \cite{ZX,PZ}. We will consider the solution $(u_\epsilon, v_\epsilon, h_\epsilon, g_\epsilon)$ of 
\begin{align}\label{cp_2}
\left\{\begin{array}{ll}
u_{\varepsilon
,t} - \Delta u_{\varepsilon} + \lambda |\nabla u_{\varepsilon}|^{\alpha} = av_{\varepsilon}^{p},\ \text{for}\ t>0,  0<r<h_{\varepsilon}(t),\\
v_{\varepsilon, t} - \Delta v_{\varepsilon} + \lambda |\nabla v_{\varepsilon}|^{\alpha} = au_{\varepsilon}^{p},\ \text{for}\ t>0,  0 < r < g_{\varepsilon}(t),\\
 u_{\varepsilon}(t,h_{\varepsilon}(t)) = v_{\varepsilon}(t,g_{\varepsilon}(t))= 0, u_{\varepsilon,r}(t,0)= v_{\varepsilon,r}(t,0) = 0,\ \text{for}\ t>0, \\  
h_{\varepsilon}^{\prime}(t) = -\mu u_{\varepsilon,r}(t,h_{\varepsilon}(t)), g_{\varepsilon}^{\prime}(t)= -\eta  v_{\varepsilon,r}(t,g_{\varepsilon}(t)),\ \text{for}\ t>0,\\
h_{\varepsilon}(0) = h_{0}^{\varepsilon}, g_{\varepsilon}(0) = g_{0}^{\varepsilon}, u_{\varepsilon}(0,r) = u_{0,\varepsilon}(r), \text{for}\ 0 \leq r \leq h_0^{\varepsilon}, v_{\varepsilon}(0, r) = v_{0,\varepsilon}(r)\ \text{for}\ 0 \leq r \leq g_0^{\varepsilon},
\end{array}\right.
\end{align}
 with $\varepsilon > 0 $ small enough $h_{0}^{\varepsilon} = h_0(1-\varepsilon)$ and $g_{0}^{\varepsilon} = g_0(1-\varepsilon)$. Furthermore, $u_{0,\varepsilon}(r)$ and $v_{0,\varepsilon}(r)$  
satisfy the followings
\begin{align*}
\left\{\begin{array}{l}
u_{0,\varepsilon}, \in C^{2}([0,h_{0}^{\varepsilon}]), v_{0,\varepsilon} \in C^{2}([0,g_{0}^{\varepsilon}]), u_{0,\varepsilon}(h_{0}^{\varepsilon})=v_{0,\varepsilon}(g_{0}^{\varepsilon})=0, \\
0 < u_{0,\varepsilon}(r) < u_0(r)\ \text{in}\ (0,h_{0}^{\varepsilon}) \ \text{and}\ 0 < v_{0,\varepsilon}(r) < v_0(r) \ \text{in} \ (0,g_{0}^{\varepsilon}),\\
u_{0,\varepsilon}(\dfrac{h_0}{h_{0}^{\varepsilon}}r) \to u_{0}(r) \ \text{in} \ C^{2}([0,h_0]) \ \text{norm as} \ \varepsilon \to 0 \ \text{and} \ v_{0,\varepsilon}(\dfrac{g_0}{g_{0}^{\varepsilon}}r) \to v_{0}(r) \ \text{in} \ C^{2}([0,g_0]) \ \text{norm as} \ \varepsilon \to 0. 
\end{array}\right.
\end{align*}
\ \ \ We obtain that $h_{\varepsilon}(t) < \overline{h}(t) $ and $g_{\varepsilon}(t) < \overline{g}(t)$ for all $t \in (0,T]$. In fact, if our assertion do not hold, there exists a $t_0 \in (0,T]$ such that $h_\varepsilon(t) <\overline{h}(t)$, $g_\varepsilon(t) <\overline{g}(t), \ \forall t \in (0,t_0)$, and either $h_\varepsilon(t_0)=\overline{h}(t_0)$ or $g_\varepsilon(t_0)=\overline{g}(t_0)$ holds. Hence 
 \begin{equation} \label{cp_3}
 h_{\varepsilon}^{\prime}(t_0) \geq \overline{h}^{\prime
 }(t_0)\ \text{or}\  g_{\varepsilon}^{\prime}(t_0) \geq \overline{g}^{\prime
 }(t_0)
\end{equation}
\ \ \ We denote  $U(t,r)= u_{\varepsilon}(t,r)-\overline{u}(t,r)$ in  $\mathcal{O}^1_{t_0}=\{ (t,r) \in \R^{2}: 0 < t \leq t_0, 0 < r < h_{\varepsilon}(t)\}$ and $V(t,r)= v_{\varepsilon}(t,r)-\overline{v}(t,r)$ in $\mathcal{O}^2_{t_0}=\{ (t,r) \in \R^{2}: 0 < t \leq t_0, 0 < r < g_\varepsilon(t)\}$ then we have 
\begin{align}\label{cp_4}
\left\{\begin{array}{l}
U_{t} - U_{rr} + A_1 U_{r} \leq B_1V,\ \text{for}\ t>0,  0<r<h_{\varepsilon}(t),\\
V_{t} - V_{rr} + A_2V_{r} \leq B_2U,\ \text{for}\ t>0,  0 < r < g_{\varepsilon}(t),\\
U(t,h_{\varepsilon}(t)) =  V(t,g_{\varepsilon}(t)) = 0,\ \text{for}\ t>0,\\
U_{r}(t,0) \geq 0, V_{r}(t,0) \geq  0, \ \text{for}\ t>0, \\  
 U(0,r) <0,\ \text{for}\ 0 \leq r \leq h_0^{\varepsilon}; V(0, r)< 0,\ \text{for}\ 0 \leq r \leq g_0^{\varepsilon},
\end{array}\right.
\end{align}
where

$
\begin{aligned}
&A_1 =A_{01}-\dfrac{N-1}{r},\ A_2 = A_{02}-\dfrac{N-1}{r},\\
&A_{01} =\lambda \alpha \int_{0}^{1}\big|sU_r+\overline{u}_r\big|^{\alpha-2}[sU_r+\overline{u}_r]ds , \  B_1 = p\int_0^1 [sV+\overline{v}]^{p-1}ds,\\
&A_{02} =\lambda \alpha \int_{0}^{1}\big|sV_r+\overline{v}_r\big|^{\alpha-2}[sV_r+\overline{v}_r]ds , \  B_2 = p\int_0^1 [sU+\overline{u}]^{p-1}ds
\end{aligned}
$

and $A_{0i}, B_{i}$ ($i=1,2$) are all nonnegative and bounded. Choose $r_0 \in (0, s_0]$ where $s_0 = \min\{h^{\varepsilon}_0,g^{\varepsilon}_0\}$ such that $U \leq 0$ and $V \leq 0$ in $[0,r_0]\times[0,t_0]$. Moreover, there exists $\mathcal{K} >0$ such that
\begin{eqnarray}
\label{cp_5a}
U_{t} -  U_{rr} \leq \mathcal{K} |U_r|+ B_1V,\ \text{for}\ t>0, r_0 <r<h_{\varepsilon}(t),\\
\label{cp_5b}
V_{t} - V_{rr} \leq  \mathcal{K} |V_r| +B_2U,\ \text{for}\ t>0, r_0 < r < g_{\varepsilon}(t).
\end{eqnarray}

 Applying Lemma \ref{lm_cp0}, we have $U \leq 0$ in $(r_0, h_\varepsilon )\times [0,t_0]$ and $V \leq 0$ in $(r_0, g_\varepsilon )\times [0,t_0]$ and then it follows from the strong maximum principle that $U < 0$ in $\mathcal{O}^1_{t_o}$ and $V < 0$ in $\mathcal{O}^2_{t_o}.$ We infer that $U_r(t_0,h_\varepsilon(t_o)) > 0$ and $V_r(t_0,g_\varepsilon(t_o))> 0$ which imply $$-\mu(u_{\varepsilon,r}(t_0,h_\varepsilon(t_o)) < -\mu(\overline{u}_{r}(t_0,\overline{h}(t_o))\ \text{and}\ -\eta v_{\varepsilon,r}(t_0,g_\varepsilon(t_o)) < -\eta \overline{v}_{r}(t_0,\overline{g}(t_o)).$$
In the other words, 
\begin{equation}\label{cp_7}
h_{\varepsilon}^{\prime}(t_0) < \overline{h}^{\prime}(t_0)\ \text{and}\ g_{\varepsilon}^{\prime}(t_0) < \overline{g}^{\prime}(t_0)\ 
\end{equation}
that is contradiction to (\ref{cp_3}). We may apply the usual comparison principle over $\mathcal{O}^1_T$ and $\mathcal{O}^2_T$ to conclude that 
\begin{equation}\label{cp_8}
u_\varepsilon(t,r) < \overline{u}(t,r)\ \text{in}\  \mathcal{O}^1_T, \  v_\varepsilon(t,r) < \overline{v}(t,r), \  \text{in}\ \mathcal{O}^2_T.
\end{equation}
Since the unique solution  $(u,v,h,g)$ of (\ref{eqs_main}) depends continuously on the initial data in (\ref{eqs_main}) and $(u_\varepsilon, v_\varepsilon,h_\varepsilon,g_\varepsilon) \to (u,v,h,g)$ as $\varepsilon \to 0,$ the desired result will be obtained by letting $\varepsilon \to 0$ in (\ref{cp_8}). The proof of Lemma \ref{lm_cp} is finished.
\end{proof}
\begin{lemma}[Monotonicity]\label{lm_monotonicity}
Let $\mathcal{D}_T = \{ (t,r) \in \R^{2}: 0 < t \leq T, 0 < r < s(t)\}$, $T >0$ and $f=f(s,y,\xi), g=g(s,y,\xi) \in C^{1}([0,\infty)\times[0,\infty)\times\R)$. Assume that $u, v \in C(\overline{\mathcal{D}_T})\cap C^{1,2}(\mathcal{D}_T)$ satisfies $u_r, v_r \in W^{1,2;2}_{loc}(\mathcal{D}_T)\cap C(\overline{\mathcal{D}_T)}$, $u_r \leq 0$, $v_r \leq 0$ on $\Gamma = (\{0\} \times [0,s_0])\cup ((0,T)\times \{0,s(t)\})$ and
\begin{eqnarray}
\label{lm_mono_1} u_t - u_{rr} - \frac{N-1}{r}u_r = f(u,v,u_r)\ \text{on}\ \mathcal{D_T},\\
\label{lm_mono_2} v_t - v_{rr} - \frac{N-1}{r}v_r = g(u,v,u_r)\ \text{on}\ \mathcal{D_T}.
\end{eqnarray}
Then $u_r \leq 0$ and $v_r \leq 0$ on $\mathcal{D}_T$.
\end{lemma}
\begin{proof} For $\lambda >0$ large enough, the function $w = u_r e^{-\lambda t}$ and $z = v_r e^{-\lambda t}$ solve the equations
 \begin{eqnarray}
\label{lm_mono_3} w_t - w_{rr} = b_1 w_r +c_1w +d_1z \ \text{on}\ \mathcal{D_T},\\
\label{lm_mono_4} z_t - z_{rr} = b_2 z_r+c_2z +d_2w \ \text{on}\ \mathcal{D_T},
\end{eqnarray}
where $b_1 = f_{\xi}(u,v,u_r)+\dfrac{N-1}{r}$, $c_1 = f_{\xi}(u,v,u_r)-\dfrac{N-1}{r^2}-\lambda \leq 0$, $d_1 = f_{y}(u,v,u_r)$, $b_2 = g_{\xi}(u,v,u_r)+\dfrac{N-1}{r}$, $c_2 = g_{\xi}(u,v,u_r)-\dfrac{N-1}{r^2}-\lambda \leq 0$  and  $d_2 = g_{y}(u,v,u_r)$.

Assume for contradiction that $2m = \sup_{\mathcal{D_T}}w > 0$ and $2n = \sup_{\mathcal{D_T}}z > 0$. Set $\mathcal{W} = w - m, \mathcal{Z} = z - n$ and choose $r_0 \in (0, s_0)$ such that $\mathcal{W} \leq m$ and $\mathcal{Z} \leq n$ in $[0,T]\times[0,r_0]$. Then, there exists $K_1, K_2 > 0$, we have, in $[0,T]\times(r_0, s(t))$

$$\begin{aligned}
\mathcal{W}_t - \mathcal{W}_{rr} &= b_1 \mathcal{W}_r +c_1(\mathcal{W}+m) +d_1(\mathcal{Z}+n)\\
& \leq K_1 |\mathcal{W}_r| +c_1\mathcal{W} +d_1\mathcal{Z}
\end{aligned}$$
and 
$$
\begin{aligned}
\mathcal{Z}_t - \mathcal{Z}_{rr} &= b_2 \mathcal{Z}_r+c_2(\mathcal{Z}+n) +d_2(\mathcal{W}+m)\\
& \leq K_2 |\mathcal{Z}_r|+c_2\mathcal{Z}+d_2\mathcal{W}.
\end{aligned}$$
Thank to Lemma \ref{lm_cp0}, we infer that $\mathcal{W} \leq 0$ and $\mathcal{Z} \leq 0$ in $[0,T]\times(r_0,s(t))$, contradicting the definition of $m, n$. If $m \leq 0$ and $n > 0$, or $n \leq 0$ and $m > 0$, there is also a contradiction. So $u_r \leq 0$ and $v_r \leq 0$ on $\mathcal{D}_T$.  The proof is comleted. 
\end{proof}
\section{Blow up in finite time}\label{section_Blowup}
In this section, we study the blowup behavior of the solution and the conditions it needs. Firstly, we give a sufficient condition under which the blowup occurs 
\\
\textbf{\textit{Proof of theorem \ref{theorem_main1}(1).}} To prove Theorem \ref{theorem_main1}(1), we combine the results of Lemma \ref{lm_ftb} and Lemma \ref{lm_ftb_2}.
\begin{lemma}\label{lm_ftb}
If \textbf{(H}) holds, $\alpha < p$ and $A$ is sufficiently large, then $T_{max} < \infty$ and there exists $C >0$ depends on $h_0, g_0, \gamma, p, q, \phi$ and $\varphi$ such that $ T_{\max} \leq  C.A^{-(p-1)}.$
\end{lemma} 
\begin{proof}
Before proving this Lemma, we observe the following functions 
\begin{center}
$W(\xi) = d^{2}(1+\dfrac{C}{2})-\dfrac{\xi^2}{2C},$ for all $\xi \in [0,dM],$
\end{center}
with $M = \sqrt{C^2+2C}$, and $d, C$ are two positive numbers to be chosen later, then it's not hard to see that the functions $W$ satisfies
\begin{align}\label{thr_fbl_1}
W(\xi) \in [d^2, d^2(1+\dfrac{C}{2})]\ \text{and}\ W'(\xi) \in [-d, 0]\ \text{for all}\ \xi \in [0, dC].
\end{align}
\begin{align}\label{thr_fbl_2}
 W(\xi) \in [0, d^2],\ \text{and}\ W'(\xi) \in [-\dfrac{dM}{C}, -d]\ \text{for all}\ \xi  \in [dC, dM].
\end{align}
\begin{align} \label{thr_fbl_3}
\begin{array}{l}
\hspace{-1.7cm} W'(\xi) = -\dfrac{\xi}{C},\ W''(\xi)= -\dfrac{1}{C},\ \text{for all}\ \xi \in (0, dM).
\end{array}
\end{align}

To prove this Lemma, we will construct a suitable subsolution  and then use comparison principle to find that $T_{\max} < \infty$ since $A$ is sufficiently large.\\
 We consider $\mathcal{O}^{i}_\varepsilon = [t^*,\frac{1}{\varepsilon})\times[0,\delta_{i}(t))$ with $\delta_{i}(t) = dM(1-\varepsilon t)^{m_{i}},\ t \in [t^*, \frac{1}{\varepsilon})$, with $i=1,2$.
\begin{center}
$\underline{u}(t,r) = \dfrac{1}{(1-\varepsilon t)^{k}}W(\dfrac{r}{(1-\varepsilon t)^{m_1}}),\ (t,r) \in \mathcal{O}^{1}_\varepsilon$,   and $\underline{v}(t,r) = \dfrac{1}{(1-\varepsilon t)^{h}}W(\dfrac{r}{(1-\varepsilon t)^{m_2}}),\ (t,r) \in \mathcal{O}^{2}_\varepsilon,$
\end{center}
where two functions $W_1, W_2$ are twice differentiable  and all positive numbers $t^*,  \varepsilon, h, k, m_1,m_2 $ shall be determined later. Moreover, 
\begin{align}\label{thr_fbl_4}
\begin{array}{rl}
L_1 &=\underline{u}_t -\underline{u}_{rr} -\dfrac{N-1}{r}\underline{u}_r + \lambda|\underline{u}_r|^{\alpha} -a\underline{v}^{p}\\
& = \bigg[\dfrac{k\varepsilon}{(1-\varepsilon t)^{k+1}}W(\xi)+ \dfrac{m_1\xi\varepsilon}{(1-\varepsilon t)^{k+1}}W'(\xi) - \dfrac{1}{(1-\varepsilon t)^{k+2m_1}}W''(\xi)\\
&\hspace{1 cm} - \dfrac{N-1}{r}\dfrac{1}{(1-\varepsilon t)^{k+m_1}}W'(\xi) + \dfrac{\lambda}{(1-\varepsilon t)^{(k+m_1)\alpha}}|W'(\xi)|^{\alpha}\bigg] - \dfrac{a W^{p}}{(1-\varepsilon t)^{hp}},
\end{array}
\end{align}
\begin{align}\label{thr_fbl_5}
\begin{array}{rl}
L_2 &=\underline{v}_t -\underline{v}_{rr} -\dfrac{N-1}{r}\underline{v}_r + \lambda|\underline{v}_r|^{\alpha} -a\underline{u}^{p}\\
& = \bigg[\dfrac{h\varepsilon}{(1-\varepsilon t)^{h+1}}W(\xi)+ \dfrac{m_2\xi\varepsilon}{(1-\varepsilon t)^{h+1}}W'(\xi) - \dfrac{1}{(1-\varepsilon t)^{h+2m_2}}W''(\xi)\\
&\hspace{1 cm} - \dfrac{N-1}{r}\dfrac{1}{(1-\varepsilon t)^{h+m_2}}W'(\xi) + \dfrac{\lambda}{(1-\varepsilon t)^{(h+m_2)\alpha}}|W'(\xi)|^{\alpha}\bigg] - \dfrac{a W^{p}}{(1-\varepsilon t)^{kp}},
\end{array}
\end{align}
 we take 
\begin{equation}\label{thr_fbl*}
k = h = \dfrac{1}{p-1}
\end{equation}
thus $k+1 =hp$ and $h+1 =kp.$ \\
Under the assumption $A$ is sufficiently large and $ d < \dfrac{s_0}{M}\ \text{and take}\ t^*\  \text{close to}\ 1/\varepsilon.$\\
With for all $ \gamma \leq 0, \beta \leq 0$, and $\alpha < p$. We denote $\gamma_{BL}= (1+dM)^{\gamma}.$ 

From (\ref{thr_fbl*}), we choose  
\begin{align}\label{thr_fbl_6a}
&0 < m_i < \min\big\{\dfrac{1}{2},\dfrac{p-\alpha}{\alpha(p-1)}\},\  C > \dfrac{h}{m_i}\ \text{and}\ 0 < \varepsilon < \dfrac{d^{2(p-1)}a_1\gamma_{BL}}{h(1+C/2)}, \text{with}\ i=1,2, 
\end{align}
which implies that $1-2m_i > 0$, $k+1 - (k+m_1)\alpha >0, h+1 - (h+m_2)\alpha >0.$ Moreover, we achieve $(r+1)^{\beta} \leq 1$ and $(r+1)^{\gamma} \geq \gamma_{BL}$.
\\
From (\ref{thr_fbl_1}) - (\ref{thr_fbl_3}), (\ref{thr_fbl_4}), (\ref{thr_fbl_5}), (\ref{thr_fbl_6a}) and condition \textbf{(H)}, we obtain\

\textit{Case 1}: $\xi \in [0,dC]$.  
 \begin{align*}
\begin{array}{ll}
L_1 &\leq \dfrac{k\varepsilon d^2(1+\dfrac{C}{2})}{(1-\varepsilon t)^{k+1}}+ \dfrac{N}{C(1-\varepsilon t)^{k+2m_1}} + \dfrac{\lambda_2 (r+1)^{\beta} d^{\alpha}}{(1-\varepsilon t)^{(k+m_1)\alpha}} - \dfrac{a_1 (r+1)^{\gamma} d^{2p}}{(1-\varepsilon t)^{hp}}\\
& \leq (1-\varepsilon t)^{-k-1}\bigg[k \varepsilon d^2(1+\dfrac{C}{2})-a_1\gamma_{BL} d^{2p}\\
&+\dfrac{N}{C}(1-\varepsilon t^*)^{1-2m_1} + \lambda_2 d^{\alpha}(1-\varepsilon t^*)^{k+1 - (k+m_1) \alpha}\bigg] \leq 0.\\
\end{array}
\end{align*}
 \begin{align*}
\begin{array}{ll}
L_2 &\leq  \dfrac{h\varepsilon d^2(1+\dfrac{C}{2})}{(1-\varepsilon t)^{h+1}} + \dfrac{N}{C(1-\varepsilon t)^{h+2m_2}} + \dfrac{\lambda_2 (r+1)^{\beta}  d^{\alpha}}{(1-\varepsilon t)^{(h+m_2)\alpha}}- \dfrac{a_1 (r+1)^{\gamma} d^{2p}}{(1-\varepsilon t)^{kp}}\\
&\leq (1-\varepsilon t)^{-h-1}\bigg[h \varepsilon d^2(1+\dfrac{C}{2})-a_1\gamma_{BL} d^{2p}\\
& +\dfrac{N}{C}(1-\varepsilon t^*)^{1-2m_2} + \lambda_2 d^{\alpha}(1-\varepsilon t^*)^{h+1 - (h+m_2)\alpha}\bigg] \leq 0.\\
\end{array}
\end{align*} 

\textit{Case 2}: $\xi \in [dC, dM].$ 
 \begin{align*}
\begin{array}{rl}
L_1 &\leq \dfrac{k\varepsilon d^2-m_1\varepsilon d^2 C}{(1-\varepsilon t)^{k+1}}+ \dfrac{N}{C(1-\varepsilon t)^{k+2m_1}} + \dfrac{\lambda_2 (r+1)^{\beta}}{(1-\varepsilon t)^{(k+m_1)\alpha}}(\dfrac{dM}{C})^{\alpha}\\
&\leq (1-\varepsilon t)^{-k-1}\bigg[\dfrac{N}{C}(1-\varepsilon t^*)^{1-2m_1}+\lambda_2(\dfrac{dM}{C})^{\alpha}(1-\varepsilon t^*)^{k+1 -(k+m_1)\alpha}\bigg] \leq 0.
\end{array}
\end{align*}
\begin{align*}
\begin{array}{rl}
L_2 &\leq \dfrac{h\varepsilon d^2-m_2\varepsilon d^2 C}{(1-\varepsilon t)^{h+1}}+ \dfrac{N}{C(1-\varepsilon t)^{h+2m_2}} + \dfrac{\lambda_2 (r+1)^{\beta}}{(1-\varepsilon t)^{(h+m_2)\alpha}}(\dfrac{dM}{C})^{\alpha}\\
&\leq (1-\varepsilon t)^{-h-1} \bigg[\dfrac{N}{C}(1-\varepsilon t^*)^{1-2m_2}+\lambda_2 (\dfrac{dM}{C})^{\alpha}(1-\varepsilon t^*)^{h+1 -(h+m_2)\alpha}\bigg] \leq 0. 
\end{array}
	\end{align*}
Consequently, one has
\begin{align*}
\begin{array}{lc}
 L_{1} \leq 0,\text{on}\ \mathcal{O}^{1}_\varepsilon, L_2 \leq 0\ \text{on}\ \mathcal{O}^{2}_\varepsilon,\ \underline{u}_r(t,0) =\underline{v}_r(t,0) = 0,\ t^* \leq t < \dfrac{1}{\varepsilon},\\
 \underline{u}(t,\delta_{1}(t)) =\underline{v}(t,\delta_{2}(t)) = 0\ \text{and}\ \delta_{1}'(t)< -\mu\underline{u}_r(t,\delta_{1}(t)), \delta_{2}'(t) <-\eta\underline{v}_r(t,\delta_{2}(t))),\ t^* \leq t < \dfrac{1}{\varepsilon},
\end{array}
\end{align*}
and we also have $\delta_{1}(0) < h_0$ and $\delta_{2}(0) < g_0$. On the other hand, for $t^*$ close to $1/\varepsilon$, it follows that 
$$ \underline{u}(t^*,r) \leq \dfrac{W(0)}{(1-\varepsilon t^*)^k} \leq A \phi(r) = u_0(r)\ \text{and}\ \underline{v} (t^*,r) \leq \dfrac{W(0)}{(1-\varepsilon t^*)^k} \leq  A \varphi(r) = v_0(r),\ r \in [0,dM],$$ 
for all $A \geq \dfrac{W(0)}{C_0(1-\varepsilon t^*)^k} $ with $C_0>0$ such that $ \min\limits_{r\in [0,dM]}\{\phi(r)\} \geq C_0$ and $\min\limits_ {r\in [0,dM]}\{\varphi(r)\} \geq C_0$. Using the comparison principle, thus
$$\underline{u}(t,r) \leq u(t-t^*,r)\  \text{on}\  \mathcal{O}^{1}_\varepsilon\ \text{and}\ \underline{v}(t,r) \leq v(t-t^*,r)\ \text{on}\ \mathcal{O}^{2}_\varepsilon.$$
Clearly, when $t^*$ as tend to $1/\varepsilon$, thus $\underline{u}(t^*,r) \to \infty$ and $\underline{v}(t^*,r) \to \infty $. So, $T_{\max}(A\phi, A\varphi) \leq 1/\varepsilon - t^*< \infty$. Moreover, this yields the sharp estimate  
\begin{equation}\label{up_Tmax}
T_{\max}(A\phi, A\varphi) \leq \dfrac{1}{\varepsilon}(\dfrac{d^2(1+C/2)}{A C_0})^{p-1},\ \text{for all}\ A \geq \dfrac{d^2(1+C/2)}{C_0(1-\varepsilon t^*)^k}. 
\end{equation}
The proof of theorem \ref{theorem_main1}(1) is finished.
\end{proof} 
Next, we shall prove some estimates, which lead to verify that all solutions  would blow up in finite time in the $L^{\infty}$ sense.  

\begin{lemma}\label{lm_ftb_1}
Assume that \textbf{(H)} holds and let $(u,v,h,g)$ be the unique positive solution to the system (\ref{eqs_main}). Suppose that $0 < u(t,r) \leq M$ in $\mathcal{O}^1_T$ and $0 < v(t,r) \leq M$ in $\mathcal{O}^2_T$ for some positive number $M$ in $\mathcal{O}^1_T = \{(t,r) \in \R^{2} : 0 \leq t < T, 0 < r < h(t)\}$ and $\mathcal{O}^2_T = \{(t,r) \in \R^{2} : 0 \leq t < T, 0 < r < g(t)\}$, for some $T \in (0, \infty)$. Then there exists a constant $K > 0$ independent of $T$ such that 
 \begin{equation}\label{lm_ftb_1_1}
 0 < h'(t) <\mu K,\ \text{and}\ 0< g'(t) \leq \eta K\ \text{for all}\ t \in [0, T).  
 \end{equation}
 Moreover,
\begin{equation}\label{lm_ftb_1_2}
 0 \leq |\nabla u| \leq H(T)\ \text{and}\  0 \leq |\nabla v| \leq H(T),  
  \end{equation}
where $H(T) = \left\{\begin{array}{ll}
\dfrac{K(T+1)}{1-TC}, &\ \text{if}\ (H)\ \text{holds with}\ \gamma <0\ \text{and}\ TC < 1,\\
K,&\ \text{if}\ (H)\ \text{holds with}\ \gamma = 0,  
\end{array}\right.$ for some  $K > 0$ independent of T and $C = a_2pM^{p-1}$. 
\end{lemma} 
\begin{proof}
The inequality $h'(t), g'(t) > 0$ is obvious by using the Hopf's lemma to system (\ref{eqs_main}). Now, we shall prove that there exists $K > 0$, independent of $T$, satisfies  (\ref{lm_ftb_1_1}). Let
$$ \mathcal{O}_{K_1} =\big\{ (t,r) \in \R^{2} : t \in (0, T), h(t) - K_1^{-1} < r < h(t)\big\}, $$
we define 
$$X(t,r) = M[2K_1(h(t) - r)- K_1^2(h(t)-r)^2]$$
with some appropriate $K_1 > 0$ to be chosen later. Direct calculations yield, for $(t,r) \in \mathcal{O}_{K_1}$, that
\begin{align*}
X_t = 2MK_1h'(t)[1-K_1(h(t)-r)];\quad\quad
X_r = 2MK_1[K_1(h(t) - r) - 1];\quad\quad X_{rr} = -2MK_1^2,
\end{align*}
and $X_t \geq 0, X_r \in (-2MK_1, 0).$ 
Assume first that $0\leq v\leq M$. Taking $K_1 \geq \max\{\dfrac{2}{h_0+1},\left(a_2\frac{M^{p-1}(h_0+1)^{\gamma}}{2^{\gamma+1}}\right)^{1/2}\} $, we have
\begin{align*}
&X_t - \Delta X + \lambda |\nabla X|^{\alpha} = 2MK_1 h'(t) [1-K_1(h(t)-r)] + 2MK_1^2 -\dfrac{N-1}{r^2}2MK_1[K_1(h(t)-r)-1]\\
& +\lambda (2MK_1)^{\alpha} [1-K_1(h(t)-r)]^{\alpha} \geq 2MK_1^2 \geq 2M (h_0+1 - K_1^{-1})^{-\gamma}K_1^2(r+1)^{\gamma} \geq av^p = u_t - \Delta u +\lambda |\nabla u|^{\alpha}.
\end{align*}
 In addition, $X(t,h(t) -K_1^{-1}) \geq u(t, h(t) - K_1^{-1}), X(t,h(t)) = u(t,h(t)) = 0$. We can choose $K_1 $ in such the way that  
$K_1 = \max \big\{\dfrac{2}{h_0+1},\frac{4\|u_0\|_{C^1([0,h_0])}}{3M},\left(a_2\frac{M^{p-1}(h_0+1)^{\gamma}}{2^{\gamma+1}}\right)^{1/2} \big\}$ then  $X(0,r) \geq u_0(r),$ for $h_0 -K_1^{-1}\leq r \leq h_0$, and we apply the maximum principle to $X - u$ over $\mathcal{O}_{K_1}$, thus 
$u(t,r) \leq X(t,r)\ \text{for}\ (t,r) \in \mathcal{O}_{K_1}$
and 
$u_r(t,h(t)) \geq X_r(t,h(t)) = -2MK_1 $ in $(0, T)$. So $0 < h'(t) \leq \mu 2MK_1$.\\
 We define the second comparison function by 
$$Y(t,r) = M[2K_2(g(t) - r)- K_2^2(g(t)-r)^2],$$
on $\mathcal{O}_{K_2} =\big\{ (t,r) \in \R^{2} : t \in (0, T), g(t) - K_2^{-1} < r <g(t)\big\}$
and choose 
 $$K_2 =  \max \big\{\dfrac{2}{g_0+1},\frac{4\|v_0\|_{C^1([0,g_0])}}{3M},\left(a_2\frac{M^{p-1}(g_0+1)^{\gamma}}{2^{\gamma+1}}\right)^{1/2} \big\}.$$
We can also prove that $v_r(t,g(t)) \geq Y_r(t,g(t)) = -2MK_2$ in $(0, T)$. So, $0 < g'(t) = -\eta v_r(t,g(t)) \leq   \eta 2MK_2$. \\

\textbf{\textit{Case 1.}} Assumption \textbf{(H)} holds with $\gamma < 0$, we can write $\lambda_{r} \leq 0$ and $|a_r| \leq a_2$.\\

Next, we shall prove (\ref{lm_ftb_1_2}), if \textbf{(H)} holds with $\gamma < 0$. Let $H(t,r) = u_r(t,r)$, $G(t,r) = v_r(t,r)$, we have 
\begin{align}\label{lm_ftb_1_4}
\left\{\begin{array}{ll}
H_t - \Delta H + \dfrac{N-1}{r^2}H = apv^{p-1}G + a_{r}v^{p}-\lambda \alpha |H|^{\alpha-2}H.H_r - \lambda_{r}|H|^{\alpha},& t>0,\  0<r<h(t),\\
G_t - \Delta G + \dfrac{N-1}{r^2}G = apu^{p-1}H + a_{r}u^{p}-\lambda \alpha |G|^{\alpha-2}G.G_r - \lambda_{r}|G|^{\alpha}, & t>0,\  0 < r < g(t),\\
H(t,0) = G(t,0)= 0,\ H(t,h(t)) < 0,\ G(t,g(t)) < 0, & t>0, \\  
h^{\prime}(t) = -\mu H(t,h(t)); g^{\prime}(t)=-\eta G(t,g(t)),& t>0,\\
 H(0,r) = u'_0(r) \leq 0,\ 0 \leq r \leq h_0\ \text{and}\ \ G(0, r) = v'_0(r) \leq 0,\ 0 \leq r \leq g_0.\\ 
\end{array}\right.
\end{align}
We write $H = e^{a_0t}W,\ G = e^{a_0t}Z$ with $a_0$ is a positive constant then, we get the following system
\begin{align}\label{lm_ftb_1_5}
\left\{\begin{array}{l}
W_t - W_{rr} = A_1 W_r + B_1 W +C_1Z + E,\ \text{for}\ t>0,\  0<r<h(t),\\
Z_t - Z_{rr} = A_2 Z_r + B_2 Z +C_2 W+ F,\ \text{for}\ t>0,\  0 < r < g(t),\\
W(t,0) = Z(t,0)= 0,\ W(t,h(t)) < 0,\ Z(t,g(t)) < 0,\ \text{for}\ t>0, \\  
W(0,r) = u'_0(r) \leq 0,\ 0 \leq r \leq h_0\ \text{and}\ \ Z(0, r) = v'_0(r) \leq 0,\ 0 \leq r \leq g_0,\\ 
\end{array}\right.
\end{align}
where 
\begin{center}
$A_1 = \dfrac{N-1}{r}-\lambda \alpha e^{a_0(\alpha-1)t}|W|^{\alpha-1}W$,  $C_1 =a p v^{p-1}$,

$A_2 = \dfrac{N-1}{r}-\lambda \alpha e^{a_0(\alpha-1)t}|Z|^{\alpha-1}Z$,  $C_2 = a p u^{p-1}$, 

$B_1 = -a_0 - \dfrac{N-1}{r^2}+ \lambda_r |W|^{\alpha-1} e^{(a_0-1)\alpha t}$, $B_2=-a_0 - \dfrac{N-1}{r^2} + \lambda_r |Z|^{\alpha-1} e^{(a_0-1) \alpha t}$,

 $E = e^{-a_0t}a_rv^p $ and $F = e^{-a_0t}a_ru^p $.
\end{center}
Assume for contradiction
that $2m =  \sup_{O^1_T}W > 0$ and $2n =  \sup_{O^2_T}Z > 0$. Set $\mathcal{W} = W - m$, $\mathcal{Z} = Z - n$ and choose $r_0 \in (0, \min\{h_0, g_0\}]$ such that $W \leq m$  and $Z \leq n$ in $[0,T] \times [0,r_0]$.  Take $a_0 > 0$ large enough such that $B m + C_1 n + E \leq 0$ and $B n+ C_2 m + F \leq 0$. Then, there exists $K_1, K_2 > 0$, in $[0,T] \times (r_0, h(t))$, we have

\begin{equation}\label{lm_ftb_1_6}
\begin{aligned}
\mathcal{W}_t - \mathcal{W}_{rr} &\leq A_1 \mathcal{W}_r + B (\mathcal{W}+m) +C_1(\mathcal{Z}+n) + E\\
& \leq K_1|\mathcal{W}_r|+ B \mathcal{W} + C_1 \mathcal{Z}
\end{aligned}
\end{equation}
similarly a.e. in $[0,T] \times (r_0, g(t))$, we have
\begin{align}\label{lm_ftb_1_7}
\mathcal{Z}_t - \mathcal{Z}_{rr} & \leq K_2 |\mathcal{Z}_r| + B \mathcal{Z} + C_2 \mathcal{W}.
\end{align}
Thank to Lemma \ref{lm_cp0}, thus  $\mathcal{W} \leq 0$ in $(r_0, h(t))\times [0,T]$ and $\mathcal{Z} \leq 0$ in $(r_0, g(t))\times [0,T]$, contradicting the definition of $m$ and $n$. If $m \leq 0$ and $n > 0$, or $n \leq 0$ and $m > 0$, there is also a contradiction. We infer that $u_r \leq 0$ in  $\mathcal{O}^1_T$ and $v_r \leq 0 $ in $\mathcal{O}^2_T$. In addition, we obtain that 

\begin{equation}\label{lm_ftb_1_8}
\sup_{\mathcal{O}^{1}_T} |u_r| \leq \sup\limits_{\Gamma_1} |u_r| + T(C\sup_{\mathcal{O}^{2}_T} |v_r| + M_0)\ \text{and}\ \sup_{\mathcal{O}^{2}_T} |v_r| \leq \sup\limits_{\Gamma_2} |u_r| + T(C\sup\limits_{\mathcal{O}^{1}_T} |v_r| + M_0),
\end{equation}

 where $\Gamma_1 = (\{0\} \times [0,h_0])\cup ((0,T)\times \{h(t)\})$, $\Gamma_2 = (\{0\} \times [0,g_0])\cup ((0,T)\times \{g(t)\})$, $C = a_2pM^{p-1}$ and $M_0 >0$ independent of $T$. Moreover, it is obvious that, 
$$\sup_{\Gamma_1} |u_r| \leq \max \big\{\sup\limits_{r \in [0,h_0]} |u_0'(r)|,\ 2MK_1\big\} \leq M_0 \ \text{and}\ \sup_{\Gamma_2} |v_r| \leq \max \big\{\sup\limits_{r \in [0,h_0]} |v_0'(r)|,\ 2MK_2 \bigg\} \leq M_0, $$

with $TC <1$ then there exists $K > 0$, independent of $T$, such that
 
 $$0 \leq |\nabla u| \leq \dfrac{K(T+1)}{1-TC}, 0 \leq |\nabla v| \leq \dfrac{K(T+1)}{1-TC},0 < h'(t) <\mu K,\ \text{and}\ 0< g'(t) \leq \eta K.$$ 
 
\textbf{\textit{Case 2.}} Assumption \textbf{(H)} holds with $\gamma = 0$, then we can write $\lambda_{r} = a_r = 0$.\\

 In this case, we shall prove (\ref{lm_ftb_1_2}) if \textbf{(H)} holds with $\gamma = 0$. From (\ref{lm_ftb_1_5}), in $\mathcal{O}_T = [0,T)\times(0, k(t))$, with $k(t) = \max\{h(t),g(t)\}$ for all $t \in [0,T)$, we get the following system in matrix form 
\begin{align}\label{lm_ftb_1_9}
\partial_t U 
-\partial_{rr} U  -  A_r(\partial_r U) + A_0 U  = \theta 
\end{align}
Where $U = \begin{pmatrix} W\\ Z \end{pmatrix},\ A_0 = \begin{pmatrix}
a + \dfrac{N-1}{r^{2}}& -pu^{p-1}\\
-pv^{p-1}& a+\dfrac{N-1}{r^{2}}
\end{pmatrix},\ \theta = \begin{pmatrix}
0\\0
\end{pmatrix},$ and
$$\ A_r = \begin{pmatrix} \lambda \alpha |W|^{\alpha -1}W -\dfrac{N-1}{r}& 0\\
0 & \lambda \alpha |Z|^{\alpha -1}Z -\dfrac{N-1}{r}
\end{pmatrix}.$$
Taking $\xi_o = \begin{pmatrix}
\xi_1^0\\ \xi_2^0
\end{pmatrix},\ \xi = \begin{pmatrix}
\xi_1\\ \xi_2
\end{pmatrix} \in \R^{2}$ and  $\langle \xi_0,\xi \rangle  = \xi_1^0.\xi_1 + \xi_2^0.\xi_2 = 0$. We obtain 
\begin{align*}
&\langle\xi_0, \xi_0\rangle + \langle A_r\xi_0, \xi \rangle + \langle A_0\xi, \xi \rangle = (\xi_1^0)^2 + \dfrac{N-1}{r^2}(\xi_1)^2+\alpha \lambda |W|^{\alpha-1}We^{a(\alpha-1)t}\xi_1^0\xi_1+ \\
& + (\xi_2^0)^2 + \dfrac{N-1}{r^2}(\xi_2)^2+\alpha \lambda |Z|^{\alpha-1}Ze^{a(\alpha-1)t}\xi_2^0\xi_2   + a(\xi_1)^{2} + a(\xi_2)^2 - p(u^{p-1}+v^{p-1})\xi_1\xi_2.
\end{align*}
If $u \equiv v$ and $a \geq pM^{p-1}$. It is easy to see that $\langle\xi_0, \xi_0\rangle + \langle A_r\xi_0, \xi \rangle + \langle A_0\xi, \xi \rangle \geq 0,$ for all $\xi_0, \xi \in \R^2$.\\
Otherwise, if there exists $a_0 >0$ such that for all $ t \in [0,T]$, we have 
$$ \sup_{[0,\ h(t)]} \dfrac{\lambda \alpha r}{2\sqrt{r^2+N-1}}|\nabla v|^{\alpha} \leq e^{a_0 (t+1)},\ \text{and}\ 
 \sup_{[0,\ g(t)]} \dfrac{\lambda \alpha r}{2\sqrt{r^2+N-1}} |\nabla v|^{\alpha} \leq e^{a_0 (t+1)}.
$$
Choosing $a \geq \max\{pM^{p-1}+1, a_0\},$
we get the following that
$$ \dfrac{(\alpha \lambda)^2 r^2 |\nabla u|^{2\alpha}}{4[(N-1)+r^2]}  \leq e^{2a(t+1)},\ 
 \dfrac{(\alpha \lambda)^2 r^2 |\nabla v|^{2\alpha}}{4[(N-1)+r^2]} \leq e^{2a(t+1)},\ \text{and}\
 pM^{p-1}+1 \leq a,
$$
thus,
$ (\alpha \lambda)^{2} |W|^{2\alpha}e^{2a(\alpha-1)t} \leq 4(\dfrac{N-1}{r^2}+1),\ 
 (\alpha \lambda)^{2} |Z|^{2\alpha}e^{2a(\alpha-1)t} \leq 4(\dfrac{N-1}{r^2}+1).$ \\
Therefore, we obtain 
\begin{align*}
&(\xi_1^0)^2 + (\dfrac{N-1}{r^2}+1)(\xi_1)^2+\alpha \lambda |W|^{\alpha-1}W e^{a(\alpha-1)t}\xi_1^0\xi_1 \geq 0,\\
&(\xi_2^0)^2 + (\dfrac{N-1}{r^2}+1)(\xi_2)^2+\alpha \lambda |Z|^{\alpha-1}Ze^{a(\alpha-1)t}\xi_2^0\xi_2 \geq 0,\\
\text{and}\ & (a-1)(\xi_1)^{2} + (a-1)(\xi_2)^2 - p(u^{p-1}+v^{p-1})\xi_1\xi_2 \geq 0. 
\end{align*}
So, $\langle\xi_0, \xi_0\rangle + \langle A_r\xi_0, \xi \rangle + \langle A_0\xi, \xi \rangle \geq 0,$ for all $\xi_0, \xi \in \R^2$. 
We apply Proposition 9.4 in \cite{VG}, then 
$$0 \leq \sqrt{|\nabla u|^{2}+|\nabla v|^{2}} \leq \max \{\sqrt{|\nabla u|^{2}+|\nabla v|^{2}} : (t,r) \in \Gamma \}$$
 with $\Gamma = (\{0\} \times [0,k_0])\cup ((0,T)\times \{k(t)\})$. In addition,
$$\max \{\sqrt{|\nabla u|^{2}+|\nabla v|^{2}} : (t,r) \in \Gamma \} \leq K: = \max \bigg\{ \big[(\sup_{r \in [0,h_0]} u_0'(r))^2+ (\sup_{r \in [0,g_0]} v_0'(r))^2 \big]^{1/2},\ 2M\sqrt{K_1^2+K_2^2} \bigg\}.$$
We deduce  $|\nabla u| \leq K$ and   $|\nabla v| \leq K$, which completes the proof. 
\end{proof}
\begin{lemma}\label{lm_ftb_2}
Suppose that \textbf{(H)} holds, $(u,v,h,g)$ be the positive solution to system (\ref{eqs_main}) with   $T_{\max} < \infty$, 
\begin{enumerate}[i)]
\item If \textbf{(H)} holds with $\gamma < 0$ and $A$ is sufficiently large, then 
 \[\limsup\limits_{t\to T_{\max}}(\|u(t,\cdot)\|_{L^{\infty}([0,t]\times[0,h(t)])} +\|v(t,\cdot)\|_{L^{\infty}([0,t]\times[0,g(t)])}) = \infty.\]
 \item If \textbf{(H)} holds with $\gamma =0$, then 
 $\limsup\limits_{t\to T_{\max}}(\|u(t,\cdot)\|_{L^{\infty}([0,t]\times[0,h(t)])} +\|v(t,\cdot)\|_{L^{\infty}([0,t]\times[0,g(t)])}) = \infty.$
\end{enumerate}

\end{lemma}
\begin{remark} We note that, since the system (\ref{eqs_main}) associated with a free boundaries, it is not directly implied that the solution blows up in finite time $T_{\max}$ when $T_{\max}<\infty$. Therefore, this lemma confirms the blow-up phenomena when $T_{\max}<\infty$ for the system (\ref{eqs_main}).

\end{remark}
\begin{proof}
i)
We use the contradiction argument that $T_{\max} < \infty$, 
$$\limsup\limits_{t\to T_{\max}}\|u(t,\cdot)\|_{L^{\infty}([0,t]\times[0,h(t)])}< \infty\ \text{and}\ \limsup\limits_{t\to T_{\max}}\|v(t,\cdot)\|_{L^{\infty}([0,t]\times[0,g(t)])} < \infty.$$
Thus, there exists $M, S>0$ such that $T_{\max} < S < \infty$ and $u(t,r) \leq M$  for all $t \in [0, T_{\max}]$, $r \in [0, h(t)]$ and $\ v(r,t) \leq M$ for all $t \in [0, T_{\max}]$, $r \in [0, g(t)]$. Note that  thank to (\ref{up_Tmax}), we can choose $A$ is suitblely large such that $T_{max} < S < \dfrac{M^{1-p}}{a_2p} < \infty$. According to Lemma \ref{lm_ftb_1}, we have $0 < h'(t) \leq \mu K$,  and $0 < g'(t) \leq \eta K$, $t \in [0, T_{\max}],$ and $0 \leq |\nabla u| \leq H(T_{\max})$, $0 \leq |\nabla v| \leq H(T_{\max})$, 
 hence $h_0 < h(t) \leq \eta KT_{\max}+h_0$ and $g_0 < g(t) \leq \eta KT_{\max}+h_0$. Moreover, $0 \leq |\nabla u| \leq H(S)$ and $0 \leq |\nabla v| \leq H(S)$ for some $K >0$ independent of $T_{\max}$. For any fixed $\gamma_0 \in (0, T_{\max}),$ by the $L^{p}$ arguments, the Sobolev imbedding and the Schauder estimates, thus there exists $M_1 > 0$, depending only on $\gamma_0, M, K$ and $S$ such that  
 \[\|u(t, \cdot)\|_{C^{2}([0,h(t)])} + \|v(t, \cdot)\|_{C^{2}([0,g(t)])} \leq M_1,\ \text{for}\ t \in [\gamma_0, T_{\max}).\]
In the same way as the proof of Proposition \ref{proposition_lc&uni}, we claim that there exists $\tau > 0$ such that the unique solution to the system (\ref{eqs_main}) with the initial time $T_{\max} - \tau$ can be extended uniquely to $T_{\max}-\tau +2\tau$, which is a contradiction to the definition of $T_{\max}$.\\
ii) The same way of argument as above with suppose \textbf{(H)} holds with $\gamma =0$, we also have 
\[\limsup\limits_{t\to T_{\max}}(\|u(t,\cdot)\|_{L^{\infty}([0,t]\times[0,h(t)])} +\|v(t,\cdot)\|_{L^{\infty}([0,t]\times[0,g(t)])}) = \infty.\] 
 The proof of Lemma \ref{lm_ftb_2} is complete.     
\end{proof}

\section{Global existence}\label{section_global}
In this section, we study global existence of positive solution to system (\ref{eqs_main}). Firstly, we prove the existence of global fast solution with suitable small initial functions.\\\
\textbf{\textit{Proof of theorem \ref{theorem_main1}(2)(i)}.}  
\begin{proof}
With $\gamma \leq 0, \beta \leq 0$, we study the following functions 
\begin{equation}\label{thr_gfs_1}
\overline{U}(t,r) = c e^{-kt}w(\dfrac{r}{\overline{s}_{1}(t)}),\ 0 \leq r < \overline{s}_{1}(t), t \geq 0\ \text{and}\ \overline{V}(t,r) = d e^{-ht}w(\dfrac{r}{\overline{s}_{2}(t)}),\ 0 \leq r < \overline{s}_{2}(t), t \geq 0,
\end{equation}
where $\overline{s}_{1}(t) = 2h_0(2 - e^{-l_{1}t})$, $\overline{s}_{2}(t) = 2g_0(2 - e^{-l_{2}t})$,  and $w(\xi) = 1- \xi^{2},\ 0 \leq \xi \leq 1$ with $c, d, k, h, l_1, l_2 >0$ will be determined later. For all $t \geq 0$ and $0 \leq r \leq \overline{s}_{1}(t)$, we obtain that 
\begin{align*}
\begin{array}{rl}
&\overline{U}_t -\Delta \overline{U} +\lambda(t,x) |\nabla \overline{U}|^{\alpha} - a(t,x)\overline{V}^{p}=\overline{U}_t -\overline{U}_{rr} -\dfrac{N-1}{r}\overline{U}_r + \lambda|\overline{U}_r|^{\alpha} -a\overline{V}^{p}\\
& \geq  c e^{-kt}[-k+k\xi^{2}+2\xi r \dfrac{\overline{s}_{1}'}{\overline{s}_{1}^{2}}+\dfrac{2}{\overline{s}_{1}^{2}}+2\xi\dfrac{N-1}{r})\dfrac{1}{\overline{s}_{1}}+ \lambda_{1}(r+1)^{\beta} c^{\alpha - 1}e^{-(\alpha-1)kt}(\dfrac{\xi}{\overline{s}})^{\alpha}\\
&-\dfrac{d^{p}}{c}a_{2}(r+1)^{\gamma}.e^{-(hp-k)t}w^{p}]\\
& \geq ce^{-kt}[-k+\dfrac{2}{\overline{s}_{1}^2}-\dfrac{d^p}{c}a_2(r+1)^{\gamma}.e^{-(hp-k)t})]\\
& \geq ce^{-kt}[-k+\dfrac{1}{8h_0^2}-\dfrac{d^p}{c}a_2.e^{-(hp-k)t}],
\end{array}
\end{align*}
similarly, we also get that $$ \overline{V}_t -\Delta \overline{V} +\lambda |\nabla \overline{V}|^{\alpha} - a\overline{U}^{p} 
 \geq de^{-ht}[-h+\dfrac{1}{8g_0^2}-\dfrac{c^p}{d}a_2.e^{-(kp-h)t}],$$
 for all  $t \geq 0$ and $0 \leq r \leq \overline{s}_{2}(t)$. Moreover,
\begin{align*} 
&\overline{s}_{1}'(t) = 2h_0 l_{1} e^{-l_{1}t},\ \overline{s}_{2}'(t) = 2g_0 l_{1} e^{-l_{2}t}; - \mu\overline{U}_r(t, \overline{s}_{1}(t)) = \dfrac{2\mu}{ \overline{s}_{1}}c e^{-kt}; -\eta \overline{V}_r(t,\overline{s}_{2}(t)) = \dfrac{2\eta}{ \overline{s}_{2}}d e^{-ht}.\\
&\overline{U}(0,r)= c\dfrac{4h_0^2-r^{2}}{4h_0^2} \geq \dfrac{3c}{4},\ \text{and}\ \overline{V}(0,r)= d\dfrac{4g_0^2-r^{2}}{4g_0^2} \geq \dfrac{3d}{4},\ 0 \leq r \leq s_0.
\end{align*}  
We will take $k = h = l_1p = l_2 = \dfrac{1}{16s_0^2}$,  $c = d \leq \min\{(\dfrac{1}{16a_2s_0^2})^{\frac{1}{p-1}}, \dfrac{h_0^2}{8ps_0^2\mu}, \dfrac{g_0^2}{8s_0^2\eta}\}$, with $s_0 = \max\{h_0,g_0\}$ then for all $t >0,$
\begin{align*}
\begin{array}{c}
\overline{U}_t -\Delta \overline{U} +\lambda |\nabla \overline{U}|^{\alpha} - a\overline{V}^{p} \geq 0,\ \text{for}\ 0 \leq r \leq \overline{s}_1(t),\\
\ \overline{V}_t -\Delta \overline{V} +\lambda |\nabla \overline{V}|^{\alpha} - a\overline{U}^{p} \geq 0,\ \text{for}\ 0 \leq r \leq \overline{s}_2(t),\\
\overline{s}_{1}'(t) \geq  - \mu\overline{U}_r(t, \overline{s}_{1}(t)\ \text{and}\ \overline{s}_{2}'(t) \geq -\eta \overline{V}_r(t,\overline{s}_{2}(t) ).
\end{array} 
\end{align*}
Moreover, $\overline{U}_r(t,0)) \leq 0,\ \overline{V}_{r}(t,0) \leq 0, \overline{U}(t,\overline{s}_1(t)) = \overline{V}(t,\overline{s}_2(t)) = 0$, $\overline{s}_1(0)=2h_0 > h_0, \overline{s}_2(0)= 2g_0 > g_0$, $\overline{U}(0,r) \geq \|u_0\|_{\infty}$ and $\overline{V}(0,r) \geq \|v_0\|_{\infty}$. Applying the comparison principle, we have

$ h(t) \leq \overline{s}_{1}(t),\ g(t) \leq \overline{s}_{2}(t)$, $u(t,r) \leq \overline{U}(t,r)$ for $t > 0$, $0 < r < h(t)$ and $v(t,r) \leq \overline{V}(t,r)$ for $t > 0$, $0 < r < g(t)$,\\
 thus
$h_{\infty} \leq \lim_{t \to +\infty}\overline{s}_1(t) = 4h_0 < +\infty$, $g_{\infty} \leq \lim_{t \to \infty}\overline{s}_2(t) = 4g_0 < \infty$  \text{and},
 $$u(t,r) \leq \overline{U}(t,r) \leq C_1 e^{-kt},\ v(t,r) \leq \overline{V}(t,r) \leq C_2 e^{-ht}.$$
The proof is complete.
\end{proof}

Next, we study the existence of a global slow solution of system (\ref{eqs_main}) by proving Theorem \ref{theorem_main1}(2)(ii). We need the following lemma.  
\begin{lemma}\label{lm_gss_1} 
Suppose that \textit{\textbf{(H)}} holds, $(u,v,h,g$) be solution of the system (\ref{eqs_main}) with $T_{max} = \infty$, $h_{\infty}< \infty$, $g_{\infty} < \infty$ and 
\begin{equation}\label{lm_gss_1_1b}
\sup_{t \geq 0} \|u(t, \cdot)\|_{L^{\infty}(0,h(t))} + \sup_{t \geq 0} \|v(t, \cdot)\|_{L^{\infty}(0,g(t))}\leq C,
\end{equation}
where $C$ is a positive constant and depends on $\|u_0\|_{C^{1}}, \|v_0\|_{C^1}$, $h_0$ and $ g_0$.
Then the solution of system (\ref{eqs_main}) satisfies  $\lim_{t \to \infty} \|u(t,\cdot)\|_{L^{\infty}(0,h(t))} =\lim_{t \to \infty} \|v(t,\cdot)\|_{L^{\infty}(0,g(t))} = 0.$ 
\end{lemma}

{ \begin{remark} In view of this lemma, we underline that the global solution of (\ref{eqs_main}) is allowed to be unbounded a priori. This result confirms that if the initial condition is small enough, the global solution must be globally bounded.
\end{remark}}

\begin{proof}
We use the contradiction argument. Assume that 
\begin{center}
$\nu:= \limsup_{t \to \infty} \|u(t,\cdot)\|_{L^{\infty}(0,h(t))} > 0 $ or  $\sigma : = \limsup_{t \to \infty} \|v(t,\cdot)\|_{L^{\infty}(0,g(t))} > 0.$
\end{center}

Case 1: With $\nu:= \limsup_{t \to \infty} \|u(t,\cdot)\|_{L^{\infty}(0,h(t))} > 0 $, then there exists a sequence $(t_n, r_n) \in (0, \infty)\times(0, h(t)) \subset (0, \infty)\times(0, h_\infty)$ such that 
$u(t_n,r_n) \geq \dfrac{\nu}{2}\ \text{for all}\ n \in \N$, where $t_n \to \infty$ as $n \to \infty$. By sequence $(r_n)$ is bounded, we have a subsequence of $(r_n)$, denote $(r_{n_k})$ that converges to $r_0 \in (0, h_\infty)$, we have $r_{n_k} \to r_0$ as $k \to \infty$.\\
 Define
$u_n(t,r) = u(t+t_n,r)$ for $t \in (-t_n,\infty),\ r \in (0, h(t+t_n))$ and $v_n(t,r) = v(t+t_n,r)$, for $t \in (-t_n,\infty),\ r \in (0, g(t+t_n))$. 

From (\ref{lm_gss_1_1b}), we know $(u_n)$ and $(v_n)$ are bounded, hence, according to the parabolic regularity it follows that $(u_n)$ there exists a subsequence $(u_{n_k})$ and $(v_n)$ there exists a subsequence $(v_{n_k})$ such that $u_{n_k} \to \overline{u}$ and $v_{n_k} \to \overline{v}$ as $k \to \infty$. Furthermore, $\overline{u}$ and $\overline{v}$ satisfy 
\[\overline{u}_t - \Delta \overline{u}  =a\overline{v}^{p} -\lambda |\nabla \overline{u}|^{\alpha},\ t \in \R,\  0 < r < h_\infty\ \text{and}\  \overline{v}_t - \Delta \overline{v} = a\overline{u}^{p} - \lambda |\nabla \overline{v}|^{\alpha},\ t \in \R,\  0 < r < g_\infty.\]
On the other hand,
 $\overline{u}(0,r_0):= \lim\limits_{k \to \infty} u_{n_k}(0,r_{n_k}) = \lim\limits_{k \to \infty}u(t_{n_k},r_{n_k}) \geq \dfrac{\nu}{2} > 0,$ and  $\overline{v}(0,r_0):= \lim\limits_{k \to \infty} v_{n_k}(0,r_{n_k}) = \lim\limits_{k \to \infty}v(t_{n_k},r_{n_k}) \geq 0,$ thus $\overline{u} > 0$ in $\R \times (0,h_\infty)$ and $\overline{v} \geq 0$ in $\R \times (0,g_\infty)$. Since 
\begin{center}
$\overline{u}(0,h_\infty) = \lim\limits_{k \to \infty} u_{n_k}(0,h(t_{n_k})) = \lim\limits_{k \to \infty}u(t_{n_k},h(t_{n_k}))=0.$
\end{center}
Using Hopf's Lemma for the above equations, we look for $\overline{u}$ at the point ($0, h_\infty$) satisfies $\overline{u}_{r}(0, h_\infty) < 0$ and from the proof of Proposition \ref{proposition_lc&uni}, we know that
\[\|u\|_{C^{\frac{1+\beta_0}{2},1+\beta_0}(0,\infty)\times(0,h(t))}+\|v\|_{C^{\frac{1+\beta_0}{2},1+\beta_0}(0,\infty)\times(0,g(t))}+\|h\|_{C^{1+\frac{\beta_0}{2}}(0,\infty)} +\|g\|_{C^{1+\frac{\beta_0}{2}}(0,\infty)}\leq C.\]
Thus $ \|h\|_{C^{1+\frac{\beta_0}{2}}(0,\infty)} \leq C$, on the other hand $h(t)$ is increasing and bounded, therefore $h'(t) \to 0$. By the Stefan condition, we have $u_r(t_n,h(t_n)) \to 0$. Moreover, 
  $\|u\|_{C^{\frac{1+\beta_0}{2},1+\beta_0}(0,\infty)} \leq C,$
it implies that $u_r(t_n,h(t_n)) = (u_n)_r(0,h(t_n)) \to \overline{u}_r(0,h_\infty)$ 
as $n \to \infty,$ which holds a contradiction with $\overline{u}_{r}(0, h_\infty) < 0.$\\
So, $\lim\limits_{t \to \infty} \|u(t,\cdot)\|_{L^{\infty}(0,h(t))} = 0.$

Case 2: With $\sigma= \limsup_{t \to \infty} \|v(t,\cdot)\|_{L^{\infty}(0,g(t))} > 0 $. Similar to above, we can deduce  a contradiction. Therefore, $\lim_{t \to \infty} \|v(t,\cdot)\|_{L^{\infty}(0,g(t))} = 0$.
The proof is complete. 
\end{proof}

\textbf{Proof of Theorem \ref{theorem_main1}(2)(ii).}   
\begin{proof}
We first denote $Aw_0 = A(\phi, \varphi)$ and the solution of (\ref{eqs_main}) by $(u(Aw_0,\cdot ) , v(Aw_0,\cdot ), h(Aw_0,\cdot ), g(Aw_0,\cdot ) )$ to emphasize the dependence of $u, v, h$ and $g$ on the initial value, so is $h_{\infty}, g_{\infty} $ and the maximal existence time $T$. We can define the following set
\begin{align*}
S = \{A > 0 : T(Aw_0) = \infty,\ h_{\infty}(Aw_0)+ g_{\infty}(Aw_0) < \infty\ \text{and}\ u(Aw_0,\cdot), v((Aw_0,\cdot)\ \text{are bounded}\}.
\end{align*}
Clearly, $S$ is not empty and bounded by Theorems \ref{theorem_main1}(1) and \ref{theorem_main1}(2)(i). Setting $A_0^* =\sup S > 0$, $U = u(A_0^*w_0, \cdot), V = v(A_0^*w_0,\cdot)$, $\delta = h(A_0^*w_0,\cdot), \sigma = g(A_0^*w_0,\cdot)$ and $\tau =  T(A_0^*w_0)$, we can get that $\tau = \infty$ and $\delta_{\infty} + \sigma_{\infty} = \infty$.\\ 
Firsly, for any fixed $t \in [0,\tau),$ by continuous dependence in Proposition \ref{proposition_cdoid}, we have $u(Aw_0,\cdot)$ converges to $U(t)$ and $v(Aw_0,\cdot)$ converges to $V(t)$ in $L^{\infty}(0, \infty)$, $h(Aw_0,t) \to \delta(t)$ and $g(Aw_0,t) \to \sigma(t)$ as $A \to A_0^*$. Since $T(Aw_0) = \infty$ for $A \in (0, A_0^*)$, we obtain that 
$\|U(t)\|_{\infty} +\|V(t)\|_{\infty} \leq C$ for all $t \in [0, \tau)$. Suppose that $\tau < \infty$, then according to lemma \ref{lm_ftb_2} $(ii)$, non-global solution must satisfy 
 $\limsup_{t \to \tau} \|u(t,\cdot)\|_{\infty} = \infty$ or $\limsup_{t \to \tau} \|v(t,\cdot)\|_{\infty} = \infty$, we infer that $\tau = \infty.$\\
Next to, clearly, we obtain that 
$\|U(t)\|_{\infty} +\|V(t)\|_{\infty} \leq C$ for all $t \geq 0$, where $C >0$ depends on $\|u_0(A^{*}_0 w_0)\|_{C^{1}}, \|v_0(A^{*}_0w_0)\|_{C^1}$, $\delta_0$ and $ \sigma_0$. Then, $\sup_{t \geq 0}\|U(t)\|_{\infty} +\sup_{t \geq 0}\|V(t)\|_{\infty} \leq C$. Now, we suppose that $\delta_{\infty} < \infty$ and $\sigma_\infty < \infty$ then $\|U(t)\|_{\infty} + \|V(t)\|_{\infty} \to 0$ as $t \to \infty$ by Lemma \ref{lm_gss_1}. On the other hand, by Theorem \ref{theorem_main1}(2)(i), there exists some large $t_0$ such that
\begin{equation}
\|U(t_0)\|_{\infty} + \|V(t_0)\|_{\infty} \leq \dfrac{3}{4} \min\{(\dfrac{1}{16a_2S_0^2})^{\frac{1}{p-1}}, \dfrac{\delta_0^2}{8pS_0^2\mu}, \dfrac{\sigma_0^2}{8S_0^2\eta}\},
\end{equation}   
where $S_0=\min\{\delta_0, \sigma_0\}$. Since $A \to A_0^*$, we have  
\begin{center}
$\|u(Aw_0,t_0)\|_{\infty} + \|v(Aw_0,t_0)\|_{\infty} \leq \dfrac{3}{4}\min\{(\dfrac{1}{16a_2(S_0^2)(Aw_0,t_0)})^{\frac{1}{p-1}}, \dfrac{h_0^2(Aw_0,t_0)}{8pS_0^2(Aw_0,t_0)\mu}, \dfrac{g^2_0(Aw_0,t_0)}{8S_0^2(Aw_0,t_0)\eta}\},$ 
\end{center}
 for $A > A_0^*$ sufficently close to $A_0^*$ by continuous dependence, thus $h_{\infty}(Aw_0)< \infty$ and $g_{\infty}(Aw_0) < \infty$ that a contradiction to the definition of $A_0^*$.
The proof is complete.
\end{proof} 
Before showing the result of theorem \ref{theorem_main1}(2)(iii), we need to prove the following lemma.  
\begin{lemma}\label{lm_gls*}
Let $x, \sigma > 0$ and $\kappa \in \R$, then

\ \  i) $x^{a}+ \dfrac{\sigma}{x^{b}} \geq \sigma ^{\frac{a}{a+b}}$, with $a, b >0$;
 
 \ \ ii) $x^{\alpha} \geq \sigma^{\frac{\alpha-s}{\alpha -1}}x^{s}-\sigma x$, and $m^{\kappa}x^{\alpha} \geq m^{\kappa\frac{s-1}{\alpha-1}}\sigma^{\frac{\alpha-s}{\alpha -1}}x^{s}-\sigma x$ with $\alpha > s > 1$ and $m >0$.     
\end{lemma}
\begin{proof}
\textit{i)} For $a, b > 0$, we consider the following function $f(x) = x^{a}+\sigma x^{-b},\ x > 0$, thus $f'(x) = a x^{a-1} - b\sigma x^{-b-1}$ and $f'(x) = 0$ since $x = \big(\dfrac{b \sigma}{a}\big)^{\frac{1}{a+b}}$. We always get that $f'(x) > 0$ for $x > \bigg(\dfrac{b \sigma}{a}\bigg)^{\frac{1}{a+b}}$ and $f'(x) < 0$ for $0 < x < \bigg(\dfrac{b \sigma}{a}\bigg)^{\frac{1}{a+b}}$. Hence, $f(x) \geq f\bigg(\bigg(\dfrac{b \sigma}{a}\bigg)^{\frac{1}{a+b}}\bigg)$,  for $x >0$, with $f\bigg(\bigg(\dfrac{b \sigma}{a}\bigg)^{\frac{1}{a+b}}\bigg) = \sigma^{\frac{a}{a+b}}\big[ \big(\dfrac{a}{b}\big)^{\frac{a}{a+b}} + \big(\dfrac{b}{a}\big)^{\frac{b}{a+b}}\big] \geq \sigma^{\frac{a}{a+b}}$ by $ \big(\dfrac{a}{b}\big)^{\frac{a}{a+b}} + \big(\dfrac{b}{a}\big)^{\frac{b}{a+b}} \geq 1.$  
\\
\textit{ii)} Applying \textit{i)} for $a = \alpha - s > 0$ and $b = s - 1$, we get that $x^{\alpha -s }+ \dfrac{\sigma}{x^{s - 1}} \geq \sigma ^{\frac{\alpha -s}{\alpha -1}}$, thus $x^{\alpha} \geq \sigma^{\frac{\alpha-s}{\alpha -1}}x^{s}-\sigma x$ and then we replace $\sigma$ by $\sigma.m^{-\kappa}$, we get $m^{\kappa}x^{\alpha} \geq m^{\kappa\frac{s-1}{\alpha-1}}\sigma^{\frac{\alpha-s}{\alpha -1}}x^{s}-\sigma x$.   
The proof is complete.
\end{proof}

\textbf{\textit{Proof of theorem \ref{theorem_main1}(2)(iii)}}    
\begin{proof} With $\gamma \leq \beta \dfrac{p - 1}{\alpha -1}$ and $\alpha \geq p$. 
We consider the following test functions, for $0 < t < T$,
\begin{align*}
\overline{u}(t,r) = C_1 e^{\gamma_0 t}[h(r)-h(\delta_1)]\,\ 0< r <\delta_2(t),  \text{and}\ \ \overline{v}(t,r) = C_2 e^{\gamma_0 t}[h(r)-h(\delta_2)],\ 0< r <\delta_2(t), 
\end{align*} 
where $h(s) = e^{-\varepsilon m s}\cosh(\varepsilon s),\ s \geq 0,$ and $m \geq 2$ is a given constant number, $\delta_{i}(t) = \gamma_{0}^{-1}C_i\varepsilon m \max\{\mu,\eta \}e^{\gamma_{0} t}, i = 1,2$ for $ 0 < t < T$, with $\varepsilon, C_1, C_2, \gamma_0$ are all positive real numbers to be determined later. For all $t \geq 0$, $s \geq 0$ and $0 < r < \delta_1(t)$, one has

\ \ \ \ \ \  $h'(s) = -\varepsilon e^{-\varepsilon ms}[m\cosh(\varepsilon s)-\sinh(\varepsilon s)] < 0,\ \delta_{i}'(t) = C_i m \varepsilon \mu_{0} e^{\gamma_{0} t}, i = 1,2,$ with $\mu_0 = \max\{ \mu, \eta\},$
\begin{align*}
&h''(s) = - \varepsilon^{2}e^{-\varepsilon m s}[-(1+m^{2})\cosh(\varepsilon s)+2m\sinh(\varepsilon s)],\\
& h'(\delta_i)\delta_i' = - C_i\varepsilon^{2} m \mu   e^{\gamma_0 t-\varepsilon m \delta_i}[m\cosh(\varepsilon \delta_i)-\sinh(\varepsilon \delta_i)],\ 0<-h'(\delta_i) < \varepsilon m \\
\text{and}\ & \gamma_0 h(\delta_i)+ h'(\delta_i)\delta_i' \leq 
 - e^{-\varepsilon m \delta_i}[(C_i\varepsilon^{2} \mu_0 (m^2/2) e^{\gamma_0 t} - \gamma_0)\cosh(\varepsilon\delta_i) + C_i\varepsilon^{2} \mu_0 e^{\gamma_0 t}m(m/2-1)\sinh(\varepsilon \delta_i)] \leq 0,
\end{align*}
by taking $C_i$ large enough such that $C_i\varepsilon^{2} \mu_0 (m^2/2) e^{\gamma_0 t} - \gamma_0 \geq 0, i= 1,2.$ Furthermore, we see that
 \begin{align*}
&[m\cosh(\varepsilon r)-\sinh(\varepsilon r)]^p \geq [(m/2)\cosh(\varepsilon r)+(m/2 - 1)\sinh(\varepsilon r)]^p \geq (m/2)^{p}\cosh^p(\varepsilon r). 
\end{align*}
Thanks to \textit{ii)} of Lemma \ref{lm_gls*}, direct computations yield 
\begin{align*}
\begin{array}{ll}
\hspace{0.4 cm} & \overline{u}_t -\Delta \overline{u} +\lambda |\nabla \overline{u}|^{\alpha} - a\overline{v}^{p} = \overline{u}_t -\overline{u}_{rr} -\dfrac{N-1}{r}\overline{u}_r + \lambda|\overline{u}_r|^{\alpha} -a\overline{v}^{p}\\
& \geq  C_1e^{\gamma_0 t}[\gamma_0 h(r) - \gamma_0 h(\delta_1) - h'(\delta_1)\delta_1'] - C_1 e^{\gamma_0 t}h''(r) -\dfrac{N-1}{r}C_1 e^{\gamma_0 t}h'(r)+\lambda_1 (r+1)^{\beta}(C_1 e^{\gamma_0 t}|h'(r)|)^{\alpha} \\
&- a_2C_2^{p}(r+1)^{\gamma}e^{p\gamma_0 t}[h(r) - h(\delta_2)]^{p}\\
&\geq C_1 e^{\gamma_0 t}\gamma_0 h(r) - C_1 e^{\gamma_0 t}h''(r)+\lambda_1 \big[(r+1)^{\beta\frac{p-1}{\alpha-1}}\sigma_1^{\frac{\alpha -p}{\alpha -1}}(C_1 e^{\gamma_0 t}|h'(r)|)^{p} - \sigma_1 C_1 e^{\gamma_0 t}|h'(r)|\big]\\
& - (r+1)^{\gamma}a_2C_2^{p}e^{p\gamma_0 t}|h(r)|^{p}
\end{array}
\end{align*}
\begin{align*}
\begin{array}{rl}
&\geq C_1e^{\gamma_0 t -\varepsilon m r}\gamma_0 \cosh(\varepsilon r)  + C_1 \varepsilon ^{2} e^{\gamma_0 t -\varepsilon m r}\big[ -(1+m^2)\cosh(\varepsilon
 r) +2m \sinh(\varepsilon r)\big]\\
&\hspace{3.3cm} - \lambda_1\sigma_1 C_1 \varepsilon e^{\gamma_0 t -\varepsilon m r}[m\cosh(\varepsilon r)-\sinh(\varepsilon r)] -(r+1)^{\gamma}a_2C_2^{p}e^{p(\gamma_0 t - \varepsilon m r)}\cosh^{p}(\varepsilon r)\\
&\hspace{6cm}+\lambda_1 (r+1)^{\beta\frac{p-1}{\alpha-1}}\sigma_1^{\frac{\alpha -p}{\alpha -1}}(C_1 \varepsilon e^{\gamma_0 t -\varepsilon m r})^{p}[m\cosh(\varepsilon r)-\sinh(\varepsilon r)]^{p}
\end{array}
\end{align*}
\begin{align*}
\begin{array}{rl}
&\geq C_1e^{\gamma_0 t -\varepsilon m r}[\gamma_0  -  \varepsilon^2(1+m^2)-\lambda_1 \sigma_1 \varepsilon m]\cosh(\varepsilon r)\\
&\hspace{1.4cm}+ \lambda_1 (r+1)^{\beta\frac{p-1}{\alpha-1}}\sigma_1^{\frac{\alpha -p}{\alpha -1}}(C_1 \varepsilon e^{\gamma_0 t -\varepsilon m r})^{p}(m/2)^{p}\cosh^p(\varepsilon r) - (r+1)^{\gamma}a_2C_2^{p}e^{p(\gamma_0 t - \varepsilon m r)}\cosh^{p}(\varepsilon r)\\
&\geq C_1e^{\gamma_0 t -\varepsilon m r}[\gamma_0  -  \varepsilon^2(1+m^2)-\lambda \sigma_1 \varepsilon m]\cosh(\varepsilon r)\\
&\hspace{1.4cm} + e^{p(\gamma_0 t - \varepsilon m r)}(r+1)^{\gamma}\big[\lambda_1 \sigma_1^{\frac{\alpha -p}{\alpha -1}}2^{-p}(mC_1 \varepsilon)^{p} - a_2C_2^{p}\big]\cosh^{p}(\varepsilon r). 
\end{array}
\end{align*}

By the same arguments as above, for all $t \geq 0$, $s \geq 0$ and $0 < r < \delta_2(t)$, we also have 
\begin{align*}
\begin{array}{rl}
&\overline{v}_t -\Delta \overline{v} +\lambda |\nabla \overline{v}|^{\alpha} - a\overline{u}^{p} = \overline{v}_t -\overline{v}_{rr} -\dfrac{N-1}{r}\overline{v}_r + \lambda|\overline{v}_r|^{\alpha} -a\overline{u}^{p}\\
&\geq C_2e^{\gamma_0 t -\varepsilon m r}[\gamma_0  -  \varepsilon^2(1+m^2)-\lambda_1 \sigma_2 \varepsilon m]\cosh(\varepsilon r)\\
&\hspace{4.6cm}+ e^{p(\gamma_0 t - \varepsilon m r)}(r+1)^{\gamma}\big[\lambda_1 \sigma_2^{\frac{\alpha -p}{\alpha -1}}2^{-p}(mC_2 \varepsilon)^{p} - a_2C_1^{p}\big]\cosh^{p}(\varepsilon r),
\end{array}
\end{align*}
where $\beta_1, \beta_2$ are two positive real numbers to be determined later.

Letting $C_1 = 2 C_2$, $\gamma_0 = \varepsilon^{2}(1+m^{2})+\lambda \max\{\sigma_1, \sigma_2\} m \varepsilon$. We consider the following cases:

$\bullet$ $\alpha > p$, we can take 
 $\sigma_1 = [(\dfrac{\lambda_1}{a_2})^{\frac{1}{p}}\varepsilon m]^{ -\frac{p(\alpha -1)}{\alpha-p}}$, $\sigma_2 = [4^{-1}(\dfrac{\lambda_1}{a_2})^{\frac{1}{p}}\varepsilon m]^{-\frac{p(\alpha -1)}{\alpha-p}}$ and $\varepsilon > 0.$
  
$\bullet$ $\alpha = p$, we can take $\sigma_1 = \sigma_2 = 1$ and $\varepsilon = (\dfrac{a_2}{\lambda_1})^{\frac{1}{p}}\dfrac{4}{m}$.\\
Then, we achieve for all $t > 0$, 
\begin{align*}
&\overline{u}_t -\Delta \overline{u} +\lambda |\nabla \overline{u}|^{\alpha} - a\overline{v}^{p} \geq 0,\ \text{for}\ 0 \leq r \leq \delta_1(t), \\
&\overline{v}_t -\Delta \overline{v} +\lambda |\nabla \overline{v}|^{\alpha} - a\overline{u}^{p} \geq 0, \ \text{for}\ 0 \leq r \leq \delta_2(t),\\
& \overline{u}(t,\delta_1(t)) = \overline{v}(t,\delta_2(t)) = 0,\ \overline{u}_r(t,0)) \leq 0,\ \overline{v}_r(t,0) \leq 0. 
\end{align*}

In addition, we will take $C_1, C_2 $ according to as above and satisfy $\overline{u}(0,r) = C_1 (h(r) -h(\delta_1(0)) \geq u_0(0,r)$ and $\overline{v}(0,r)= C_2 (h(r) -h(\delta_2(0)) \geq  v_0(0,r) $, with $\delta_i(0) = C_i\varepsilon m \mu_0 \gamma_0^{-1}, i=1,2$ and $C_1, C_2$ large enough to
\begin{center}
$\delta_1'(t) = C_1\varepsilon m \mu_0 e^{\gamma_0 t} \geq \mu C_1e^{\gamma_0 t}(-h'(\delta_1)) =  - \mu\overline{u}_r(t, \delta_1(t))$ ; 

$\delta_2'(t) = C_2\varepsilon m \mu_0 e^{\gamma_0 t}\geq C_2\eta e^{\gamma_0 t}(-h'(\delta_2))=-\eta \overline{v}_r(t,\delta(t)))$ and $\delta_1(0) > h_0, \delta_2(0) > g_0$. 
\end{center}
Therefore, with $T^*$ denotes the maximal exstence time of the supersolution $(\overline{u}(t,r), \overline{v}(t,r), \delta_1(t), \delta_2(t))$, the comparison principle yields for $0 < t < T^*$, 
$h(t) \leq \delta_1(t)$; $g(t) \leq \delta_2(t)$; 
\begin{align*}
&0 \leq u(t,r) \leq \overline{u}(t,r)\ \text{for}\ 0 \leq r \leq h(t);\\
 &0 \leq v(t,r) \leq \overline{v}(t,r)\ \text{for}\ 0 \leq r \leq g(t),
\end{align*}
 Clearly, $T^* = \infty$, thus $T_{\max} = \infty$. The proof of theorem \ref{theorem_main1}(2)(iii) is finished.   
\end{proof}

\textbf{Conflict of interest statement :} The authors have no conflicts of interest to declare that are relevant to the content of this article.



\begin{thebibliography}{99}
\bibitem {AP} A. Attouchi, P. Souplet, Gradient blow-up rates and sharp gradient estimates for diffusive Hamilton-Jacobi equations,
Calc. Var. Partial Differential Equations 59 (2020), no. 5, Paper No. 153, 28 pp.

{ \bibitem{AZ1} B. Abdelhedi,H. Zaag, Construction of a blow-up solution for a perturbed nonlinear heat equation with a gradient and a non-local term. J. Differential Equations 272 (2021), 1–45.}



{ \bibitem{AZ2} B. Abdelhedi,H. Zaag, Single point blow-up and final profile for a perturbed nonlinear heat equation with a gradient and a non-local term. Discrete Contin. Dyn. Syst. Ser. S 14 (2021), no. 8, 2607–2623.}

\bibitem {AAS} A. Al-Elaiw, S. Tayachi, Different asymptotic behavior of global solutions for a parabolic system with nonlinear gradient terms, J. Math. Anal. Appl. 387 (2012), no. 2, 970–992.

\bibitem {DMJ} D. Andreucci, M. A. Herrero and J. J. L. Velizquez, Liouville theorems and blow up behaviour in semilinear reaction diffusion systems, Ann. Inst. H. Poincaré Anal. Non Linéaire. 14 (1997), no. 1, 1–53. 

\bibitem {DE}  D. Andreucci, E. Di Benedetto, On the cauchy problem and initial traces for a class of evolution equations with strongly nonlinear sources, Ann. Scuola Norm. Sup. Pisa Cl. Sci. (4) 18 (1991), no. 3, 363–441.   


\bibitem {KKS} K. Ishige, K. Nakagawa, Salani, Paolo Spatial concavity of solutions to parabolic systems, Ann. Sc. Norm. Super. Pisa Cl. Sci. (5) 20 (2020), no. 1, 291–313.



\bibitem {SBP} S. Benachour, B. Roynette, P. Vallois, Asymptotic estimates of solutions of $u_t - \dfrac{1}{2} \Delta u = - |\nabla u|$ in $\R_{+} \times \R^d, d \geq 2$, J. Funct. Anal. 144 (1997), no. 2, 301-324.     

\bibitem {MF} M. Chipot, F.B. Weissler, Some blowup results for a nonlinear parabolic equation with a gradient term, SIAM J. Math. Anal. 20 (1989), no. 4, 886–907. 

\bibitem {RM} R. Castillo, M. Loayza, On the critical exponent for some semilinear reaction–diffusion systems on general domains, J. Math. Anal. Appl. 428 (2015), no. 2, 1117–1134. 

\bibitem {LCE} L.C. Evans, A strong maximum principle for reaction-diffusion systems and a weak convergence scheme for reflected stochastic differential equations, (2010) Massachusetts Institute of Technology. 

\bibitem {MM} M. Escobedo, M.A. Herrero, Boundedness and blow up for a semilinear reaction – diffusion system, J. Differ. Equ. 89 (1991), no. 1, 176–202. 



\bibitem {MM1} M. Escobedo, M.A. Herrero, A uniqueness result for a semilinear reaction–diffusion system, Proc. Am. Math. Soc. 112 (1991), no. 1, 175–185. 


\bibitem {MFJL} M. Fila and J. Lankeit, Continuation beyond interior gradient blow-up in a semilinear parabolic equation, Math. Ann. 377 (2020), 317–333.

\bibitem {YKH} Y. Fujishima, K. Ishige, H. Maekawa, Blow-up set of type I blowing up solutions for nonlinear parabolic systems, Math. Ann. 369 (2017), no. 3-4, 1491–1525.

\bibitem {HPD}  H. Ghidouche, P. Souplet, D.  Tarzia, Decay of global solutions, stability and blowup for a reaction – diffusion problem with free boundary, Proc. Amer. Math. Soc. 129 (2001), no. 3, 781–792. 

\bibitem {MS} M. Ghergua, S. D. Taliaferro, Pointwise bounds and blow-up for systems of semilinear parabolic inequalities and nonlinear heat potential estimates, J. Functional Analysis, 272 (2017), no. 4, 1301–1339. 


\bibitem {YG} Y. Giga, A bound for global solutions of semilinear heat equations, Comm. Math. Phys. 103 (1986), no. 3, 415–421.  

\bibitem {YR} Y. Giga, R. Kohn, Nondegeneracy of blowup for semilinear heat equations, Comm. Pure Appl. Math. 42 (1989), no. 6, 845-884.



{ \bibitem{TZ}S. Tayachi; H. Zaag, Existence of a stable blow-up profile for the nonlinear heat equation with a critical power nonlinear gradient term. Trans. Amer. Math. Soc. 371 (2019), no. 8, 5899–5972.}

\bibitem {TVH1} T.E. Ghoul, V.T. Nguyen, H. Zaag, Construction and stability of blowup solutions for a non-variational semilinear parabolic system, Ann. Inst. H. Poincaré Anal. Non Linéaire 35 (2018), no. 6, 1577–1630.

{ \bibitem {TVH2} T.E. Ghoul, V.T. Nguyen, H. Zaag, Blowup solutions for a nonlinear heat equation involving a critical power nonlinear gradient term, J. Differential Equations 263 (2017), no. 8, 4517–4564.}

\bibitem {KP} K. Ishige, P. Salani, Parabolic power concavity and parabolic boundary value problems, Math. Ann. 358 (2014), no. 3-4, 1091–1117.

\bibitem {KT} K. Ishige, T. Kawakami, Reﬁned asymptotic proﬁles for a semilinear heat equation, Math. Ann. 1 (2012), no. 353, 161–192.   
    
\bibitem {LSV} O. A. Ladyzenskaja, V. A. Solonnikov, N.N. Ural’Ceva, Linear and quasilinear equations of parabolic type. (Russian) Translated from the Russian by S. Smith. Translations of Mathematical Monographs, Vol. 23 American Mathematical Society, Providence, R.I. (1968) xi+648 pp. 

\bibitem{ALV} H. A. Levine, A  Fujita  type global existence-global nonexistence theorem for a weakly  coupled system of reaction-diffusion  equations, Z. Angew. Math. Phys. 42 (1991), no. 3, 408–430.  

\bibitem{SHCWX} S. Lian, H. Yuan, C. Cao,  W. Gao, X. Xu, On the Cauchy problem for the evolution $p$-Laplacian equations with gradient term and source, J. Differential Equations. 235 (2007), no. 2, 544–585.    

\bibitem {FH} F. Merle, H. Zaag, A Liouville theorem for vector-valued nonlinear heat equations and applications, Math. Ann. 316  (2000), no.1,  103–137.



\bibitem {NFJ} N. Mizoguchi, F. Quirós, J.L.  Vázquez, Multiple blow-up for a porous medium equation with reaction, Math. Ann.  350 (2011), no. 4, 801–827.

\bibitem {PM} P. Meier, On the critical exponent for reaction–diffusion equations, Arch. Ration. Mech. Anal. 109 (1990), 63–71.


\bibitem {VG} V. Mazya, G. Kresin, Maximum Principles and Sharp Constants for Solutions of Elliptic and Parabolic Systems, (2012) American Mathematical Society.


\bibitem {RGP1} R. G. Pinsky, Existence and nonexistence of global solutions for $u_t = \Delta u + a(x)u^p$ in $\R^d$, J. Differential Equations. 133 (1997), no. 1, 152-177. 

\bibitem {RGP2} R. G. Pinsky, The behavior of the lifespan for solutions to $u_t = \Delta u + a (x) u^p$ in $\R^d$, J. Differential Equations. 147 (1998), no. 1, 30-57.    



\bibitem {RGP4} R. G. Pinsky, Decay of mass for the equation $u_t = \Delta u-a(x)u^p|\nabla u|^q$, J. Differential Equations. 165(2000), no. 1, 1-23.

\bibitem {PQPS} P. Quitnner, P. Souplet, Superlinear parabolic problems. Blow-up, global existence and steady states, (2007) Birkhauser Advanced Texts, ISBN: 978-3-7643-8441-8.

\bibitem {RW} J. D. Rossi N. Wolanski, Blow-up vs. global existence for a semilinear reaction-diffusion system in a bounded domain, Communications in Partial Differential Equations, 20 (1995), no. 11-12, 1991-2004. 




\bibitem {PS} P. Souplet, S. Tayachi, and F. B. Weissler, Exact self-similar blow-up of solutions of a semilinear parabolic equation with a nonlinear gradient term, Indiana Univ. Math. J. 45 (1996), no. 3, 655–682. 

\bibitem {TA} A. S. Tersenov, The preventive effect of the convection and of the diffusion in the blow-up phenomenon for parabolic equations, Ann. Inst. H. Poincaré Anal. Non Linéaire 21 (2004), no. 4, 533–541.

\bibitem {SH} S. Tayachi and H. Zaag, Existence of a stable blow-up profile for the nonlinear heat equation with a critical power nonlinear gradient term.  Trans. Amer. Math. Soc. 371 (2019), no. 8, 5899–5972.

  
\bibitem {MY} M. Wang, Y. Zhao, A semilinear parabolic system with a free boundary, Z. Angew. Math. Phys. 66 (2015), no. 6, 3309–3332.



\bibitem {ZX} Z. Zhang, X. Zhang, Asymptotic behavior of solutions for a free boundary problem with a nonlinear gradient absorption, Calc. Var. Partial Differential Equations. 58 (2019), no. 1, 32:31. 

\bibitem {PZ} P. Zhou, Z.G. Lin, Global existence and blowup of a nonlocal problem in space with free boundary, J. Funct. Anal. 262 (2012), 3409–3429.      

\end{thebibliography}
\end{document}